\documentclass[a4paper,11pt,]{elsarticle}

\usepackage{amsmath,amssymb,amsthm}
\usepackage{mathtools,bm,mathrsfs}
\usepackage{graphicx,booktabs,tabularx,array}
\usepackage{placeins}
\usepackage[normalem]{ulem}
\usepackage{xcolor}
\usepackage{hyperref}
\usepackage{comment}
\usepackage{epstopdf}

\graphicspath{{figures/}}

\hypersetup{
  colorlinks=true,
  linkcolor=blue,
  citecolor=red,
  urlcolor=blue
}

\usepackage[a4paper,margin=2.5cm]{geometry}

\newtheorem{thm}{Theorem}[section]
\newtheorem{lem}[thm]{Lemma}

\newtheorem{rem}[thm]{Remark}

\numberwithin{equation}{section}

\newcommand*{\vv}[1]{\vec{\mkern0mu#1}}
\newcommand{\bR}{{\mathbb R}}
\newcommand{\Gt}{{\Gamma(t)}}
\newcommand{\Gm}{{\Gamma^m}}
\newcommand{\nabs}{\nabla_{\!s}}
\newcommand{\rd}{\;{\rm d}}
\newcommand{\id}{{\rm id}}
\newcommand{\dd}[1]{\frac{{\rm d}}{{\rm d}#1}}
\newcommand{\ddt}{\dd{t}}
\newcommand{\ttau}{\tau}
\newcommand{\nn}{\nonumber}
\newcommand{\ipd}[1]{\bigl(#1\bigr)}
\newcommand{\norm}[1]{\lVert#1\rVert}
\newcommand{\mA}{\mathscr A}
\newcommand{\mV}{\mathscr V}
\newcommand{\mVV}{\mathscr{\vv V}\rule{0pt}{0pt}}
\newcommand{\Vh}{\mathbb V^h}
\newcommand{\pq}{\partial_{\varkappa} q}
\newcommand{\pqg}{\partial_{\varkappa_\gamma} q}
\newcommand{\npq}{\nabla_{\vv p} q}

\begin{document}
\begin{frontmatter}
\title{A geometric-IEQ approach for curvature-driven flows}
\author[reg]{Harald Garcke}
\ead{harald.garcke@ur.de}
\author[tue]{Yifei Li}
\ead{yifei.li@mnf.uni-tuebingen.de}
\author[trento]{Robert N\"urnberg}
\ead{robert.nurnberg@unitn.it}
\author[ustc]{Quan Zhao}
\ead{quanzhao@ustc.edu.cn}

\address[reg]{Fakult\"at f\"ur Mathematik, Universit\"at Regensburg,
93040 Regensburg, Germany}
\address[tue]{Mathematisches Institut, Universit\"at T\"ubingen,
Auf der Morgenstelle 10, 72076 T\"ubingen, Germany}
\address[trento]{Dipartimento di Matematica, Universit\`a di Trento,
38123 Trento, Italy}
\address[ustc]{School of Mathematical Sciences, University of Science and
Technology of China, 230026 Hefei, Anhui, China}

\begin{abstract}
We introduce a geometric invariant energy quadratization (IEQ) formulation
for the $L^2$-gradient flow of a class of curvature-dependent energies of
closed hypersurfaces, where the energy density depends on the mean curvature
and the unit normal.  A
square-root auxiliary variable converts the geometric energy into a quadratic
functional.  Combining its evolution equation with a skew-symmetric treatment
of surface transport yields a weak formulation whose energy law depends only
on the normal velocity and is therefore independent of the chosen tangential
motion.  We construct monolithic nonlinear parametric finite element
methods by coupling a common IEQ discretization with one of two geometric formulations that accommodate either the Barrett-Garcke-N\"urnberg (BGN) tangential velocity or the minimal-deformation-rate (MDR) tangential velocity.  Both fully discrete approximations satisfy an
unconditional decay estimate for the discrete quadratic energy. We further extend the construction to anisotropic Willmore flow, and to $H^{-1}$-gradient flows. A variety of numerical experiments are presented to demonstrate the robustness and accuracy of the introduced methods. 
\end{abstract}


\end{frontmatter}

\section{Introduction}
\label{sec:intro}

The motion of interfaces driven by a law for the normal velocity which involves curvature
quantities is fundamental in differential geometry and arises in image
processing, materials science and models for biological membranes.  Two prototypes are mean
curvature flow, the $L^2$-gradient flow of surface area \cite{Huisken84flow}, and the Willmore flow,
the $L^2$-gradient flow of the surface integral of the squared mean curvature
\cite{Willmore93}.  For either prototype, the $L^2$-gradient-flow structure
yields a quantitative dissipation law: the time derivative of the
driving energy equals the negative squared $L^2$-norm of the normal velocity. This dissipation law is the precise geometric structure that a numerical approximation should inherit. The same variational principle admits two natural generalizations. Replacing the isotropic energy with an orientation-dependent energy leads to anisotropic mean curvature or Willmore flows~\cite{Hoffman1972vector, PPR15anisoWillmore}, whereas replacing the $L^2$-metric with an $H^{-1}$-metric yields volume conserving geometric evolutions, with surface diffusion as a classical example~\cite{Mullins57theory,Cahn94overview}. 
The $H^{-1}$-gradient flow of the Willmore energy has previously been
considered in several works~\cite{Gurtin2002interface, Hausser05discrete,  BURGER2007183,MR2500806}. Throughout this paper,
we refer to this flow as the ``$H^{-1}$ Willmore flow''.
The obtained evolution equations range from
second- to sixth-order nonlinear geometric PDEs, and their stiffness and
strong coupling of geometry and curvature make the construction of robust
energy-stable discretizations a central numerical challenge.

Parametric finite element methods (PFEMs) have a long history for approximating curvature-driven interface evolution, see \cite{Dziuk90, BMN05, BGN08parametric, Barrett20, KovacsLL21, Hu22evolving, BZ21SPFEM, BLani23, Gao26dualMDR}. Dziuk introduced a foundational parametric approximation for mean curvature flow \cite{Dziuk90}. The BGN approach provides a variational framework that extends naturally to a broad class of geometric flows \cite{BGN07variational,BGN08parametric, BGN08ani3d,BGN08willmore,Barrett20}. It relies on a weak curvature equation which can induce a beneficial tangential motion to yield asymptotically equidistributed curves and conformal polyhedral surfaces in three space dimensions. These properties highlight the importance of a suitable tangential velocity for maintaining mesh quality. Alternative treatments include formulations based on the minimal-deformation-rate (MDR) principle \cite{Hu22evolving, Gao26dualMDR}, which control mesh deformation and alleviate small-time-step mesh instabilities in BGN schemes, as well as harmonic-map/DeTurck reparametrizations and uniformization strategies \cite{DeTurck17,Duan24new}.

A separate challenge is to preserve the gradient-flow structure at the
discrete level.  Variational BGN formulations yield discrete energy
dissipation for a range of second- and fourth-order geometric flows
\cite{BGN07,BGN07variational,BGN08parametric}; for surface diffusion, structure-preserving
variants additionally preserve the enclosed volume
\cite{Jiang21, BZ21SPFEM}.  Further constructions have also been developed for
anisotropic geometric flows~\cite{BGN08Ani,BLani23,BGN08ani3d}.  For Willmore flow, parametric methods with semidiscrete energy stability were developed in early work \cite{Dziuk08,BGN08willmore} and, more recently, in \cite{Zhao26Onsager} with a relaxed MDR tangential velocity. 
More recently, Bao and Li \cite{BaoL25} developed an unconditionally energy-stable fully discrete PFEM for closed planar curves by incorporating an evolution equation for the mean curvature. The scheme in \cite{GNZ25willmore} combines this mechanism with a BGN tangential motion. The normal--tangential velocity-splitting technique of \cite{GNZ26splitting} separates the curvature-driven normal velocity, which carries the gradient-flow dissipation, from a freely chosen tangential velocity, thereby yielding an energy-stable fully discrete method for Willmore flow of curves and surfaces with arbitrary tangential motion. This development raises a natural question: can the same separation of energy dissipation and mesh motion be achieved systematically for gradient flows of more general surface energies?

These energy-stable methods for Willmore flow \cite{BaoL25, GNZ25willmore, GNZ26splitting} share two features. First, they rely on an evolution equation for the mean curvature. Second, the energy stability holds for a modified energy, in which the mean curvature is determined by this evolution equation rather than from the discrete surface. Both features are the hallmarks of the invariant energy quadratization (IEQ) methodology for classical dissipative systems \cite{Yang16IEQ,ZhaoWY17IEQ}, where an auxiliary variable is evolved and the energy stability holds for a modified energy given by this auxiliary variable. On a fixed domain, the IEQ approach introduces an auxiliary variable of the form $r=\sqrt{F(u)+C}$ for a nonlinear energy density $F$ that is bounded from below. The energy becomes quadratic in $r$, while a chain-rule evolution equation couples $r$ to the original phase variable. Evaluating the nonlinear coefficient explicitly, but using the same coefficient in the two coupled equations, produces the cancellations needed for linear, unconditionally energy-stable time discretizations. The approach has been widely applied and further developed for phase-field and other dissipative systems \cite{YangZWS17IEQ,ZhaoEtAl18EQ,GongZW20EQ,LiuLi19IEQ}.

To see this correspondence clearly, we write $g(u)=\sqrt{2F(u) + C}$. We compare the Willmore flow and the IEQ approach in Table~\ref{tab:willmore_ieq}. In both cases, the energy is quadratic in the auxiliary variable $\lambda$, namely $\lambda=\varkappa$ for Willmore flow and $\lambda=g(u)$ for IEQ. The difference is that \(g(u)\) is always non-negative whereas \(\varkappa\) is allowed to change sign. Moreover, the evolution equation for $\lambda$ is linear in the rate of change of the state variable, namely the normal velocity $\mV$ of $\Gamma$ for Willmore flow and $\partial_t u$ for IEQ. The numerical method dissipates not the original energy $E$ but the modified energy $\frac{1}{2}\int_\Gamma \lambda^2$ for Willmore flow and $\frac{1}{2}\int \lambda^2$ for IEQ. The essential difference is that the modified energy in IEQ is integrated over the fixed domain, whereas for Willmore flow it is integrated over an evolving surface.

\begin{table}[htbp]
  \centering
  \caption{Correspondence between the Willmore flow and the IEQ approach. Here $\Gamma$ is the evolving surface, $\mVV$ is the full velocity of $\Gamma$, $\vv \nu$ and $\varkappa$ are the unit normal and mean curvature of $\Gamma$, $\mV = \mVV\cdot\vv\nu$ is the normal velocity, $\Delta_s$ is the surface Laplacian, and $\partial_t^\square$ is the normal time derivative.}
  \label{tab:willmore_ieq}
  \begin{tabular}{lcc}
  \toprule
   & Willmore flow & IEQ \\
  \midrule
  State variable & $\Gamma$ & $u$ \\[1ex]
  Energy functional & $E(\Gamma)=\dfrac{1}{2}\displaystyle\int_{\Gamma} \varkappa^2$ & $E(u)=\dfrac{1}{2}\displaystyle\int g^2(u)$ \\[2ex]
  Evolution of state variable & $\mVV, \, \mV = \mVV\cdot\vv\nu$ & $\partial_t u$ \\[1ex]
  Auxiliary variable & $\lambda=\varkappa$ & $\lambda=g(u)$ \\[1ex]
  Evolution equation & $\partial_t^\square\lambda = (\Delta_s + \lvert\nabs\vv\nu\rvert^2)\,\mV$ & $\partial_t\lambda=(g'(u))\,\partial_t u$ \\
  \bottomrule
  \end{tabular}
  \end{table}

Motivated by this connection, we aim to extend the IEQ principle from
fixed-domain dissipative systems to curvature-driven geometric gradient flows on evolving
surfaces with a general surface energy.  Our main focus is the $L^2$-gradient flow of a class of
surface energies determined by the mean curvature $\varkappa$ and the unit
normal $\vv\nu$.  We develop a geometric-IEQ formulation and the corresponding
energy-stable parametric finite element approximations, whose dissipation
structure is independent of the tangential reparametrization.  We then extend
the framework in two directions.  First, we consider energies depending on
more general geometric quantities, such as the
anisotropic mean curvature $\varkappa_\gamma$, leading to 
anisotropic Willmore flow as a representative example.
Second, we treat $H^{-1}$-gradient flows, including surface-diffusion-type
evolutions, for which the resulting schemes dissipate the energy while also
preserving the enclosed volume. Throughout, it will be desirable that the method adapts to any tangential motion to ensure good mesh quality of $\Gamma$.

The main contributions of this work are as follows.
\begin{itemize}
\item We introduce a geometric-IEQ formulation for the $L^2$-gradient flow of a class of surface energies whose density depends on
  $\varkappa$ and $\vv\nu$.  The formulation uses the square root of the energy density as an auxiliary variable and
  separates its dissipative normal motion from tangential
reparametrizations. 
\item We discretize the geometric-IEQ formulation to construct nonlinear monolithic geometric--IEQ parametric finite element approximations with either BGN or MDR tangential velocities, where an unconditional
  discrete energy-decay estimate can be established for both schemes.
\item We extend the framework beyond energies that depend only on
$\varkappa$ and $\vec\nu$. For anisotropic Willmore flow we allow
a dependence on the anisotropic mean curvature $\varkappa_\gamma$, while for $H^{-1}$-gradient
  flows the resulting schemes dissipate the energy and exactly conserve the
  enclosed volume.
\end{itemize}

The remainder of the paper is organized as follows.  In
Section~\ref{sec:mf} we introduce the geometric notation and derive the
geometric--IEQ formulation for the $L^2$-gradient flow of an energy of the form $\frac12\int_\Gamma f(\varkappa,\vv\nu)$. Section~\ref{sec:pfem} presents the fully
discrete PFEMs and their energy stability estimate.  The extension to
anisotropic Willmore flow is developed in Section~\ref{sec:aniso}, followed
by the $H^{-1}$-gradient flow in Section~\ref{sec:h-1}.
Numerical results are presented in Section~\ref{sec:num-general},
including convergence tests, the validation of energy dissipation and
mass conservation, and morphological evolutions of closed curves and surfaces.

\section{Geometric--IEQ formulation for the \texorpdfstring{$L^2$}{L2}-gradient flow}
\label{sec:mf}

\subsection{Geometric setting}
\label{subsec:geometry}

Let $T>0$ and let $\{\Gt\}_{t\in[0,T]}$ be a family of smooth, closed and
oriented hypersurfaces in $\bR^d$, $d\in\{2,3\}$.  We assume that
\begin{equation}
  \vv x:\Upsilon\times[0,T]\to\bR^d
\end{equation}
is a smooth parametrization over a fixed reference manifold $\Upsilon$ and
write
\begin{equation}
  \mVV(\vv x(\vv\rho,t),t)=\partial_t\vv x(\vv\rho,t)\quad(\vv\rho, t)\in\Upsilon\times[0,T]
\end{equation}
for the induced full velocity under the parameterization.  If $\vv\nu$ is the chosen unit normal,
then the normal velocity is given by
\begin{equation}
  \mV=\mVV\cdot\vv\nu\qquad\mbox{on}\quad\Gt.
\end{equation}

For a scalar function $g$ on the evolving surface,  we define the material and
normal time derivatives on $\Gt$ through
\begin{equation}
  (\partial_t^\circ g)(\vv x(\vv\rho,t),t)
  =\ddt g(\vv x(\vv\rho,t),t)\quad(\vv\rho, t)\in\Upsilon\times[0,T],
  \qquad
  \partial_t^\square g
  =\partial_t^\circ g-\mVV\cdot\nabs g,\label{eq:timeder}
\end{equation}
see \cite[Definition~28]{Barrett20}.
We have the curvature identity
\begin{equation}
  \varkappa\,\vv\nu=\Delta_s\vec\id
  \qquad\text{on }\Gt,
  \label{eq:curvature-convention}
\end{equation}
where $\vec\id$ is the identity function in $\bR^d$. We use the sign convention such that the mean curvature of the sphere of radius $r$ is $-\frac{d-1}{r}$ with outward normal.

We recall the following standard identities; see, for example,
\cite[Lemmas~37,38,39]{Barrett20}.
\begin{lem}
For a smooth evolving hypersurface,
\begin{subequations}\label{eq:normal-evolution}
\begin{alignat}{2}
  \partial_t^\square\vv\nu&=-\nabs\mV\qquad &&\mbox{on}\quad\Gt,\label{eq:normal-evolution-nu}\\
  \partial_t^\square\varkappa
  &=\Delta_s\mV+\mV\lvert\nabs\vv\nu\rvert^2\qquad &&\mbox{on}\quad\Gt.
  \label{eq:normal-evolution-kappa}
\end{alignat}
\end{subequations}
Moreover,
\begin{subequations}\label{eq:material-evolution}
\begin{alignat}{2}
  \partial_t^\circ\vv\nu&=-(\nabs\mVV)^T\vv\nu\qquad &&\mbox{on}\quad\Gt,\\
  \partial_t^\circ\varkappa
  &=\Delta_s\mV+\mV\lvert\nabs\vv\nu\rvert^2
    +\mVV\cdot\nabs\varkappa\qquad &&\mbox{on}\quad\Gt.
\end{alignat}
\end{subequations}
\end{lem}

For later use, if $\vv a$ is a vector field on a closed hypersurface, then the surface divergence theorem reads
\begin{equation}
  \int_\Gamma \nabs\cdot\vv a
  +\varkappa\,\vv a\cdot\vv\nu\rd\Gamma=0.
  \label{eq:surface-divergence}
\end{equation}
The integration-by-parts formula on a closed hypersurface is
\begin{equation}
  \int_\Gamma\vv a\cdot\nabs g\rd\Gamma
  =-\int_\Gamma g\left(\nabs\cdot\vv a
    +\varkappa\vv a\cdot\vv\nu\right)\rd\Gamma.
\end{equation}
The surface transport theorem reads \cite[Theorem 32]{Barrett20}
\begin{align}
  \ddt\int_{\Gt}g\rd\Gamma
  &=\int_{\Gt}\left(\partial_t^\circ g + g\,\nabs\cdot\mathscr{\vv V}\right)\rd\Gamma =\int_\Gt\left(\partial_t^\circ g + g\nabs\vec\id:\nabs\mathscr{\vv V}\right)\rd\Gamma\nn\\
  &=\int_{\Gt}\left(\partial_t^\square g-g\varkappa\mV\right)\rd\Gamma,
  \label{eq:transport}
\end{align}
where we have used the identity $\nabs\cdot\mVV=\nabs\vec\id:\nabs\mVV$. 

\subsection{Geometric--IEQ reformulation of the energy}
\label{subsec:ieq}

We consider the geometric energy
\begin{equation}
  E(\Gamma)=\frac12\int_\Gamma f(\varkappa,\vv\nu)\rd\Gamma,
  \label{eq:general-energy}
\end{equation}
where $f:\bR\times(\bR^d\setminus\{\vv0\})\to\bR_{\geq0}$ is smooth and,
for every fixed $\varkappa$, positively homogeneous of degree one:
\begin{equation}
  f(\varkappa,\alpha\vv p)=\alpha f(\varkappa,\vv p)
  \qquad\forall\alpha>0,\quad\vv p\in\bR^d\setminus\{\vv 0\}.
  \label{eq:f-homogeneous}
\end{equation}
Assume that a smooth, possibly signed, factor $q$ is available such that
\begin{equation}
  q(\varkappa,\vv p)^2=f(\varkappa,\vv p),
  \qquad
  q(\varkappa,\alpha\vv p)=\alpha^{1/2}q(\varkappa,\vv p),
  \quad \alpha>0,
  \label{eq:q-factor}
\end{equation}
and define
\begin{equation}
  \lambda=q(\varkappa,\vv\nu).
  \label{eq:lambda-definition}
\end{equation}
Then
\begin{equation}
  E(\Gamma)=\frac12\int_\Gamma\lambda^2\rd\Gamma.
  \label{eq:quadratic-energy}
\end{equation}
We denote the derivatives of $q$ by $\pq(\cdot,\cdot)$ and
$\npq(\cdot,\cdot)$.  By the one-half homogeneity in
\eqref{eq:q-factor}, Euler's identity yields
\begin{equation*}
  \nabla_{\vv p}q(\varkappa,\vv p)\cdot\vv p
  =\frac12q(\varkappa,\vv p).
\end{equation*}

Using \eqref{eq:transport}, \eqref{eq:normal-evolution} and
\eqref{eq:surface-divergence}, we obtain
\begin{align}
  \ddt E(\Gt) &=\int_{\Gt} \left(\lambda\partial_t^\square q - \frac{1}{2}\lambda^2\varkappa\,\mathscr{V}\right)\rd\Gamma\nn\\
  &=\int_{\Gt}\lambda\left[
      q_\varkappa(\varkappa,\vec\nu)\left(\Delta_s\mV
      +\mV\lvert\nabs\vv\nu\rvert^2\right)
      -\npq(\varkappa,\vv\nu)\cdot\nabs\mV\right]\rd\Gamma
     -\frac12\int_{\Gt}\lambda^2\varkappa\mV\rd\Gamma
      \nn\\
  &=\int_{\Gt}\mV\left[
      \Delta_s(\lambda \pq(\varkappa,\vec\nu))
      +\lambda \pq(\varkappa,\vec\nu)\lvert\nabs\vv\nu\rvert^2
      +\nabs\cdot(\lambda\npq(\varkappa,\vv\nu))
      \right]\rd\Gamma,
  \label{eq:energy-variation}
\end{align}
where in the second equality we use the chain rule to compute 
\begin{align}
\partial_t^\square q & = q_\varkappa(\varkappa,\vec\nu)\partial_t^\square\varkappa+ \npq(\varkappa,\vv\nu)\cdot\partial_t^\square\vec\nu\nn\\
&= q_\varkappa(\varkappa,\vec\nu)\left(\Delta_s\mV
      +\mV\lvert\nabs\vv\nu\rvert^2\right)
      -\npq(\varkappa,\vv\nu)\cdot\nabs\mV.
      \label{eq:dtq}
\end{align}
On noting \eqref{eq:energy-variation}, a weighted $L^2$-gradient flow of the energy is given by
\begin{equation}
  \mV=-\beta(\vec\nu)\,\left[\Delta_s(\lambda \pq(\varkappa,\vec\nu))
  +\lambda \pq(\varkappa,\vec\nu)\lvert\nabs\vv\nu\rvert^2
  +\nabs\cdot(\lambda\,\npq(\varkappa,\vv\nu))\right]\quad\mbox{on}\quad\Gt,
  \label{eq:ieq-normal-velocity}
\end{equation}
where $\beta(\vv\nu)$ is the orientation-dependent mobility.  We assume $\beta\in C^2(\bR^d\setminus\{\vv0\})\cap C(\bR^d,\bR_{\geq0})$
and $\beta>0$ on $\bR^d\setminus\{\vv0\}$. 
For $\beta\equiv 1$ we obtain the standard $L^2$-gradient flow
of $E(\Gamma(t))$.

\subsection{Geometric--IEQ weak formulation}
\label{subsec:l2-weak}

It follows from \eqref{eq:timeder} and \eqref{eq:dtq} that 
\begin{align}
\partial_t^\circ\lambda &= \partial_t^\square\lambda + \mVV\cdot\nabs\lambda=\partial_t^\square q + \mVV\cdot\nabs\lambda\nn\\
&=\pq(\varkappa,\vec\nu)\left(\Delta_s\mV
    +\mV\lvert\nabs\vv\nu\rvert^2\right)
   -\npq(\varkappa,\vv\nu)\cdot\nabs\mV
   +\mVV\cdot\nabs\lambda, \label{eq:ieq-lambda-strong}
\end{align}
which we will use together with \eqref{eq:ieq-normal-velocity} to approximate the gradient flow. Let us next derive a geometric--IEQ weak formulation for it. 
Denote by $\ipd{\cdot,\cdot}_\Gamma$ the $L^2$-inner product on
$\Gamma(t)$.  Using integration by parts, it is not difficult to obtain that
\begin{equation}
  \ipd{\mVV\cdot\nabs\lambda,~\chi}_\Gamma
  +\frac12\ipd{\nabs\vec\id:\nabs\mVV,~\lambda\chi}_\Gamma
  =\frac12\mA_\Gamma(\mVV,\lambda,\chi)
  -\frac12\ipd{\mV,\varkappa\lambda\chi}_\Gamma
  \label{eq:transport-skew}
\end{equation}
where we introduced the skew-symmetric form
\begin{equation}
  \mA_\Gamma(\mVV,\eta,\chi)
  =(\mVV\cdot\nabs\eta,\chi)_\Gamma
  -(\mVV\cdot\nabs\chi,\eta)_\Gamma
\qquad\forall\eta, \chi\in H^1(\Gt),
  \label{eq:skew-form}
\end{equation}
satisfying $\mA_\Gamma(\mVV,\eta,\eta)=0$.  

Now the geometric--IEQ weak formulation is given as follows.  Given $\Gamma(0)$ and $\lambda(\cdot,0)$, for each $t\in(0,T]$, we find
\[
  \bigl(\mV(\cdot,t),\lambda(\cdot,t)\bigr)
  \in H^1(\Gt)\times H^1(\Gt)
\]
such that
\begin{subequations}\label{eq:general-weak}
\begin{align}
  &\ipd{[\beta(\vec\nu)]^{-1}\,\mV,~\varphi}_\Gamma
  =\ipd{\nabs(\lambda \pq(\varkappa,\vec\nu)),~\nabs\varphi}_\Gamma
   -\ipd{\lambda \,\pq(\varkappa,\vec\nu)\lvert\nabs\vv\nu\rvert^2,~\varphi}_\Gamma\nn\\
   &\qquad 
   +\frac12\ipd{\lambda^2\varkappa,~\varphi}_\Gamma
   +\ipd{\lambda\nabla_{\vv p}q(\varkappa,\vv\nu),~\nabs\varphi}_\Gamma\qquad\forall\varphi\in H^1(\Gt);
  \label{eq:general-weak-v}\\
  &\ipd{\partial_t^\circ\lambda,~\chi}_\Gamma
  +\frac12\ipd{\nabs\vec\id:\nabs\mVV,~\lambda\chi}_\Gamma
   =-\ipd{\nabs(\pq(\varkappa,\vec\nu)\,\chi),~\nabs\mV}_\Gamma +\ipd{\pq(\varkappa,\vec\nu)\,\chi\mV,
    \lvert\nabs\vv\nu\rvert^2}_\Gamma\nn\\
    &\qquad 
   -\frac12\ipd{\mV,~\varkappa\lambda\chi}_\Gamma
   -\ipd{\nabla_{\vv p}q(\varkappa,\vv\nu)\cdot\nabs\mV,~\chi}_\Gamma +\frac12\mA_\Gamma(\mVV,\lambda,\chi)\qquad\forall\chi\in H^1(\Gt).
  \label{eq:general-weak-lambda}
\end{align}
\end{subequations}

The system \eqref{eq:general-weak} determines the normal velocity $\mV$ and
the auxiliary variable $\lambda$ once the full velocity $\mVV$ and the mean
curvature $\varkappa$ are available.  To obtain a complete coupled
formulation, we supplement these geometric--IEQ equations with a normal-velocity
constraint and a geometric equation that recovers the curvature from the
parametrization. 

\begin{rem}\label{rem:onlykappa}
When the energy density depends only on the curvature, i.e., $f(\varkappa,\vec\nu)=\mathcal{F}(\varkappa)$, the homogeneity \eqref{eq:f-homogeneous} gives $f(\varkappa,\vv p)=\mathcal F(\varkappa)|\vv p|$, and the construction in \eqref{eq:q-factor} gives
\[
q(\varkappa,\vec p)=h(\varkappa)\sqrt{|\vec p|},
\qquad h(\varkappa)^2=\mathcal{F}(\varkappa).
\]
Hence $\pq(\varkappa,\vv\nu)=h'(\varkappa)$ and $\npq(\varkappa,\vec\nu)=\tfrac12 h(\varkappa)\vec\nu$. Since $\vv\nu\cdot\nabs\varphi=0$ for every scalar $\varphi$, the terms involving $\npq$ in \eqref{eq:general-weak-v} and \eqref{eq:general-weak-lambda} vanish, and \eqref{eq:general-weak} simplifies to
\begin{subequations}\label{eq:onlykappa-weak}
\begin{align}
  &\ipd{[\beta(\vec\nu)]^{-1}\,\mV,~\varphi}_\Gamma
  =\ipd{\nabs(\lambda h'(\varkappa)),~\nabs\varphi}_\Gamma
   -\ipd{\lambda h'(\varkappa)\lvert\nabs\vv\nu\rvert^2,~\varphi}_\Gamma
   +\frac12\ipd{\lambda^2\varkappa,~\varphi}_\Gamma,\\
  &\ipd{\partial_t^\circ\lambda,~\chi}_\Gamma
  +\frac12\ipd{\nabs\vec\id:\nabs\mVV,~\lambda\chi}_\Gamma
   =-\ipd{\nabs(h'(\varkappa)\chi),~\nabs\mV}_\Gamma
   +\ipd{h'(\varkappa)\chi\mV,~\lvert\nabs\vv\nu\rvert^2}_\Gamma\nn\\
   &\qquad -\frac12\ipd{\mV,~\varkappa\lambda\chi}_\Gamma
   +\frac12\mA_\Gamma(\mVV,\lambda,\chi).
\end{align}
\end{subequations}
For example, for the generalized Willmore energy 
$\frac{1}{2}\int_{\Gamma} (\varkappa - \bar\varkappa)^2$, 
with $\bar{\varkappa}\in\mathbb R$ the spontaneous curvature, we have
$\mathcal{F}(\varkappa)=(\varkappa-\bar{\varkappa})^2$, so that 
$h(\varkappa)=\varkappa - \bar{\varkappa}$, $h'(\varkappa)=1$ and $\lambda=q(\varkappa,\vec\nu)=(\varkappa-\bar{\varkappa})$. Then \eqref{eq:onlykappa-weak} reduces to the weak formulation considered in \cite{GNZ26splitting}.
\end{rem}

\subsection{Geometric formulations and tangentially independent energy law}
\label{subsec:l2-geometry-energy}

The geometric--IEQ system \eqref{eq:general-weak} needs to be
coupled with equations for $\mVV$ and $\varkappa$. We consider two
possibilities.
The first choice is the \textbf{BGN geometric formulation}
\cite{Barrett20}: find
\[
  (\mVV,\varkappa)\in[H^1(\Gt)]^d\times H^1(\Gt)
\]
such that
\begin{subequations}\label{eqn:general-weakBGN}
\begin{align}
\label{eq:general-weakBGN1}
&\ipd{\mathscr{\vv V}\cdot\vv\nu,~\xi}_\Gamma-\ipd{\mathscr{V},~\xi}_\Gamma=0\qquad\forall\xi\in H^1(\Gt);\\
&\ipd{\varkappa\,\vv\nu,~\vec\eta}_{\Gt}+\ipd{\nabs\vec\id,~\nabs\vec\eta}_\Gamma=0\qquad\forall\vv\eta\in [H^1(\Gt)]^d.
\label{eq:general-weakBGN2}
\end{align}
\end{subequations}
Alternatively, we use the \textbf{MDR geometric
formulation} \cite{Hu22evolving} by finding
\[
  (\mVV,\varkappa,\mu)
  \in[H^1(\Gt)]^d\times H^1(\Gt)\times L^2(\Gt)
\]
such that
\begin{subequations}\label{eqn:general-weakMDR}
\begin{align}
\label{eq:general-weakMDR1}
&\ipd{\mathscr{\vv V}\cdot\vv\nu,~\xi}_\Gamma-\ipd{\mathscr{V},~\xi}_\Gamma=0\qquad\forall\xi\in H^1(\Gt);\\
\label{eq:general-weakMDR2}
&\ipd{\varkappa,~g}_\Gamma+\ipd{\nabs\vv\id,\nabs(g\,\vec\nu)}_\Gamma=0\qquad\forall g\in H^1(\Gt);\\
 &\ipd{\mu\,\vec\nu,~\vec\eta}_\Gamma
 + \ipd{\nabs\mathscr{\vv V},~\nabs\vec\eta}_{\Gamma}=0
 \qquad\forall\vec\eta\in[H^1(\Gt)]^d.
\label{eq:general-weakMDR3}
\end{align}
\end{subequations}

The second equation in \eqref{eqn:general-weakMDR} is obtained by testing the curvature identity \eqref{eq:curvature-convention} with a scalar function.  The last
equation is the Euler--Lagrange equation for minimizing
\begin{equation}
  \frac12\int_{\Gt}|\nabs\mVV|^2\rd\Gamma
  \label{eq:mdr-functional}
\end{equation}
subject to the normal-velocity constraint
\eqref{eq:general-weakMDR1}, with $\mu$ acting as a Lagrange multiplier.
At the continuous level, \eqref{eqn:general-weakBGN} leaves the tangential
component of $\mVV$ undetermined; its parametric finite element discretization
selects the characteristic BGN tangential redistribution.  In contrast, the
MDR equations determine a tangential velocity through the minimization
principle \eqref{eq:mdr-functional}.  Consequently, \eqref{eq:general-weak} together with
\eqref{eqn:general-weakBGN} defines the {\bf gIEQ--BGN} system, whereas
\eqref{eq:general-weak} together with \eqref{eqn:general-weakMDR} defines the
{\bf gIEQ--MDR} system.

\begin{thm}[continuous energy law]
\label{thm:continuous-energy}
Let $(\Gt,\lambda,\mV,\mVV,\varkappa)$ be a sufficiently smooth solution of
the {\bf gIEQ--BGN} system, or let
$(\Gt,\lambda,\mV,\mVV,\varkappa,\mu)$ be a sufficiently smooth solution of
the {\bf gIEQ--MDR} system.  Then
\begin{equation}
  \ddt\left(\frac12\norm{\lambda}_{L^2(\Gt)}^2\right)
  =-\norm{[\beta(\vec\nu)]^{-\frac{1}{2}}\mV}_{L^2(\Gt)}^2\leq0.
  \label{eq:continuous-energy-law}
\end{equation}
\end{thm}
\begin{proof}
Choose $\varphi=\mV$ in \eqref{eq:general-weak-v} and $\chi=\lambda$ in
\eqref{eq:general-weak-lambda}.  All geometric coupling terms cancel and
$\mA_\Gamma(\mVV,\lambda,\lambda)=0$.  Hence
\begin{equation*}
  \norm{[\beta(\vec\nu)]^{-\frac{1}{2}}\mV}_{L^2(\Gt)}^2
  +(\partial_t^\circ\lambda,\lambda)_{\Gt}
  +\frac12(\nabs\vec\id:\nabs\mVV,\lambda^2)_{\Gt}=0.
\end{equation*}
The last two terms equal
$\frac{\rm d}{{\rm d}t}\frac12\int_{\Gt}\lambda^2\rd\Gamma$ by
\eqref{eq:transport}.  Hence the energy law holds for both geometric
formulations.
\end{proof}

\begin{rem}
It is noteworthy that the stability estimate in \eqref{eq:continuous-energy-law} relies solely on \eqref{eq:general-weak}; neither \eqref{eqn:general-weakBGN} nor \eqref{eqn:general-weakMDR} is involved in its derivation. This freedom allows the tangential motion to be chosen independently, without affecting the underlying geometric gradient flow.
\end{rem}

\section{Parametric finite element approximations of the \texorpdfstring{$L^2$}{L2}-gradient flow}
\label{sec:pfem}

We now discretize the geometric--IEQ formulation of Section~\ref{sec:mf}.
The gIEQ--BGN and gIEQ--MDR schemes share the same discrete equations for
the normal velocity and the auxiliary variable.  We give their common
energy stability proof here.

\subsection{Discrete surfaces and geometric quantities}

Let $0=t_0<t_1<\cdots<t_M=T$ be a uniform partition with time step
$\ttau$.  At time $t_m$, let
$\Gm=\bigcup_{j=1}^J\overline{\sigma}_j^m$ be a closed polyhedral
hypersurface with labelled vertices $\{\vv q_k^m\}_{k=1}^N$ and mesh
$\mathcal T^m=\{\sigma_j^m\}_{j=1}^J$, where the $\sigma_j^m$ are mutually
disjoint $(d-1)$-simplices in $\bR^d$.  We define
\begin{equation}
  \Vh(\Gm)=\{\chi\in C(\Gm):\chi|_\sigma\text{ is affine for every simplex }
  \sigma\in \mathcal{T}^m\}.
  \label{eq:discrete-space}
\end{equation}

We introduce the $L^2$-inner product and the mass-lumped inner product over the current polyhedral surface $\Gamma^m$ by
\begin{subequations}
\label{eq:surface-inner-products}
\begin{align}
  \ipd{f,~g}_{\Gamma^m}
  &:={}
  \sum_{j=1}^J
  \int_{\sigma_j^m}f\,g\,\rd\Gamma,
  \label{eq:l2mass}
  \\
  \ipd{f,~g}_{\Gamma^m}^h
  &:={}
  \frac{1}{d}\sum_{j=1}^J|\sigma_j^m|
  \sum_{k=0}^{d-1}
  \lim_{\sigma_j^m\ni\vec x\to\vec q^m_{j_k}}
  f(\vec x)\, g(\vec x),\label{eq:mass-lumping}
\end{align}
\end{subequations}
for piecewise continuous functions $f,g$ with possible jumps across the edges of $\sigma_j^m$,  where $\{\vec q^m_{j_k}\}_{k=0}^{d-1}$ are the vertices of
$\sigma_j^m$, and
$|\sigma_j^m|$ is its
surface area.  Both inner products in \eqref{eq:surface-inner-products}
extend naturally to vector- and tensor-valued functions.

Given $\Gm$, we write $\vv\nu^m$ for the element unit normal of $\Gm$, and the new parametrization is represented on the old
surface by
\newcommand{\vvX}{\vv X\rule{0pt}{0pt}}
\newcommand{\vvF}{\vv F\rule{0pt}{0pt}}
\newcommand{\vvN}{\vv N\rule{0pt}{0pt}}
\begin{equation}
  \vvX^{m+1}=\vv\id+\ttau\mVV^{m+1},
  \qquad \Gamma^{m+1}=\vvX^{m+1}(\Gm).
  \label{eq:nonlinear-surface-map}
\end{equation}
On each element of $\mathcal T^m$, set
\begin{equation}
  \vvF^{m+1}
  =\nabs\vvX^{m+1}+\vv\nu^m\otimes\vv\nu^m,
  \qquad
  \vvN^{m+1}
  =\operatorname{cof}(\vvF^{m+1})\,\vv\nu^m,
  \qquad
  \mathcal J^{m+1}=\lvert\vvN^{m+1}\rvert.
  \label{eq:nonlinear-normal}
\end{equation}
Here $\operatorname{cof}(\vvF^{m+1})$ denotes the cofactor matrix,
and
$\mathcal J^{m+1} = |\vvN^{m+1}| = \frac{|\sigma^{m+1}_j|}{|\sigma^m_j|}$ is the elementwise ratio between the measures of
$\sigma_j^{m+1}$ and $\sigma_j^m$.

We define the vertex normal $\vv v^m, \hat{\vv v}^{m+1}\in[\Vh(\Gm)]^d$ for $\Gm$ and $\Gamma^{m+1}$ by
\begin{equation}
  (\vv\omega^m,\vv\psi^h)^h_{\Gm}
  =(\vv\nu^m,\vv\psi^h)_{\Gm},
  \quad\forall\vv\psi^h\in[\Vh(\Gm)]^d,
  \quad
  \vv v^m=I_h^m\left(\frac{\vv\omega^m}{\lvert\vv\omega^m\rvert}\right),
  \label{eq:vertex-normal}
\end{equation}
\begin{equation}
  (\mathcal J^{m+1}\vv\omega^{m+1},\vv\psi^h)^h_{\Gm}
  -(\vvN^{m+1},\vv\psi^h)_{\Gm}=0,
  \quad\forall\vv\psi^h\in[\Vh(\Gm)]^d, 
  \quad
  \hat{\vv v}^{m+1}=I_h^{m}\left(\frac{\vv\omega^{m+1}}{\lvert\vv\omega^{m+1}\rvert}\right).
  \label{eq:implicit-vertex-normal}
\end{equation}
with $I_h^m$ the nodal interpolation operator onto $\Vh(\Gm)$.

We evaluate the curvature argument of $\pq$ implicitly, while retaining
the normalized vertex normal of the known surface:
\begin{equation}
  \widehat Q_\varkappa^{m+1}
  =I_h^m \pq(\varkappa^{m+1},\vv v^m),
  \qquad
  \mathcal W^m=\lvert\nabs\vv v^m\rvert^2\approx\lvert\nabs\vv\nu\rvert^2.
  \label{eq:explicit-qkappa-coefficients}
\end{equation}
Thus only the normal argument and the curvature reaction coefficient
$\mathcal W^m$ are frozen.  The interpolant
$\widehat Q_\varkappa^{m+1}\in\Vh(\Gm)$ can be represented by an additional
scalar finite element unknown with the defining equation
\begin{equation}
  (\widehat Q_\varkappa^{m+1},\zeta^h)^h_{\Gm}
  =(\pq(\varkappa^{m+1},\vv v^m),\zeta^h)^h_{\Gm}
  \qquad\forall\zeta^h\in\Vh(\Gm).
  \label{eq:implicit-qkappa-closure}
\end{equation}
In particular, $\nabs\widehat Q_\varkappa^{m+1}$ is the gradient of this
nodal interpolant, not the derivative of a differently interpolated
coefficient.  The same unknown is used in both geometric--IEQ equations.  

We also treat $\npq$ implicitly, now using both the new curvature and
the new vertex normal $\hat{\vv v}^{m+1}$ of
\eqref{eq:implicit-vertex-normal}:
\newcommand{\widehatvvQ}{\widehat{\vv Q}\rule{0pt}{0pt}}
\begin{equation}
  \widehatvvQ_{\vv p}^{m+1}
  =I_h^m\,\npq\bigl(\varkappa^{m+1},\hat{\vv v}^{m+1}\bigr).
  \label{eq:semiimplicit-qp}
\end{equation}
Although $\widehatvvQ_{\vv p}^{m+1}$ depends nonlinearly on the new geometry, the same coefficient is used in both geometric--IEQ equations.  This paired use is the only property needed for the discrete energy cancellation.

\subsection{Geometric--IEQ schemes}
\label{subsec:nonlinear-pfem}

The geometric--IEQ equations are identical for both the BGN and MDR variants. For each $m\geq 0$, we first seek 
\[
  (\mV^{m+1}, \lambda^{m+1})\in\Vh(\Gm)\times\Vh(\Gm)
\]
such that 
\begin{subequations}\label{eq:nonlinear-ieq-common}
\begin{align}
  &(\beta(\vv\nu^m)^{-1}\mV^{m+1},\varphi^h)^h_{\Gm}
  -(\widehat Q_\varkappa^{m+1}\nabs\lambda^{m+1}
    +\lambda^{m+1}\nabs\widehat Q_\varkappa^{m+1},\nabs\varphi^h)_{\Gm}
  \nn\\
  &\quad +(\lambda^{m+1}\widehat Q_\varkappa^{m+1},
    \mathcal W^m\varphi^h)_{\Gm}
  -\frac12((\lambda^{m+1})^2,
    \varkappa^{m+1}\varphi^h)^h_{\Gm}
  \nn\\
  &\quad -(\lambda^{m+1}\widehatvvQ_{\vv p}^{m+1},
    \nabs\varphi^h)_{\Gm}=0,
  \label{eq:nonlinear-ieq-v}\\
  &\left(\frac{\lambda^{m+1}-\lambda^m}{\ttau},
    \chi^h\right)_{\Gm}
  +\frac12\left(
    \nabs\vvX^{m+1}:\nabs\mVV^{m+1},\lambda^{m+1}\chi^h\right)_{\Gm}
  \nn\\
  &\quad -\frac12\mA_{\Gm}
    (\mVV^{m+1},\lambda^{m+1},\chi^h)
  \nn\\
  &\quad +(\widehat Q_\varkappa^{m+1}\nabs\chi^h
    +\chi^h\nabs\widehat Q_\varkappa^{m+1},\nabs\mV^{m+1})_{\Gm}
  -(\widehat Q_\varkappa^{m+1}\chi^h,
    \mathcal W^m\mV^{m+1})_{\Gm}
  \nn\\
  &\quad +\frac12(\mV^{m+1},
    \lambda^{m+1}\varkappa^{m+1}\chi^h)^h_{\Gm}
  +(\widehatvvQ_{\vv p}^{m+1},
    \nabs\mV^{m+1}\chi^h)_{\Gm}=0,
  \label{eq:nonlinear-ieq-lambda}
\end{align}
\end{subequations}
for every $(\varphi^h,\chi^h)\in\Vh(\Gm)^2$. Equation~\eqref{eq:nonlinear-ieq-v} approximates the normal-velocity
equation. In \eqref{eq:nonlinear-ieq-lambda}, we 
follow \cite[(4.22)]{GNZ26splitting} for the approximations for the first two terms, while the skew form uses the new parametrization velocity.

The coefficient $\widehat Q_\varkappa^{m+1}$ is coupled through
\eqref{eq:implicit-qkappa-closure}, with its normal argument explicit,
whereas $\widehatvvQ_{\vv p}^{m+1}$ uses the implicit vertex normal
through \eqref{eq:semiimplicit-qp}.  Identical coefficients
and quadrature rules are used in every cancelling pair.  In particular, the
two curvature reaction terms are both mass lumped.  The products involving
$\lambda^{m+1}$, $\varkappa^{m+1}$ and $\mVV^{m+1}$ make
\eqref{eq:nonlinear-ieq-common} a genuinely nonlinear system.

\begin{rem}\label{rem:discrete-onlykappa}
For the curvature-dependent energies described in Remark~\ref{rem:onlykappa}, the discrete terms involving
$\npq$ in \eqref{eq:nonlinear-ieq-common}  do not generally vanish when
$\widehatvvQ_{\vv p}^{m+1}$ defined by \eqref{eq:semiimplicit-qp} was used.
Using the explicit elementwise normal $\vec\nu^m$ in
\eqref{eq:semiimplicit-qp} would preserve the orthogonality property
that holds at the continuous level. However, we retain the implicit
approximation to maintain a unified scheme.
\end{rem}

Next, $\mVV^{m+1}$ and $\varkappa^{m+1}$ are determined as follows, with either the BGN tangential velocity or the MDR tangential velocity.  For gIEQ--BGN, the common equations \eqref{eq:nonlinear-ieq-common} are coupled with the vector curvature
equation.  We seek
\[(\mVV^{m+1}, \varkappa^{m+1})\in[\Vh(\Gm)]^d\times\Vh(\Gm)\]
such that, for every
$(\xi^h,\vv\eta^h)\in\Vh(\Gm)\times[\Vh(\Gm)]^d$,
\begin{subequations}\label{eq:nonlinear-bgn-geometry}
\begin{align}
  (\mVV^{m+1}\cdot\vv\nu^m,\xi^h)^h_{\Gm}
  -(\mV^{m+1},\xi^h)^h_{\Gm}&=0,
  \label{eq:nonlinear-bgn-normal}\\
  (\varkappa^{m+1}\vv\nu^m,\vv\eta^h)^h_{\Gm}
  +(\nabs\vvX^{m+1},\nabs\vv\eta^h)_{\Gm}&=0.
  \label{eq:nonlinear-bgn-curvature}
\end{align}
\end{subequations}
The first equation identifies the scalar normal velocity with the normal
component of the parametrization velocity.  The second is the PFEM
discretization of $\varkappa\vv\nu=\Delta_s\id$ on the new parametrization.
It simultaneously recovers $\varkappa^{m+1}$ and generates the characteristic
BGN tangential redistribution.  Both occurrences of the normal constraint
use the same mass-lumped pairing.

For gIEQ--MDR, we introduce a multiplier $\mu^{m+1}\in\Vh(\Gm)$ and seek
\[(\varkappa^{m+1},\mVV^{m+1},\mu^{m+1}) \in\Vh(\Gm)\times[\Vh(\Gm)]^d\times\Vh(\Gm)\]
 such that, for every
$(\xi^h,g^h,\vv\eta^h)\in\Vh(\Gm)\times\Vh(\Gm)
\times[\Vh(\Gm)]^d$,
\begin{subequations}\label{eq:nonlinear-mdr-geometry}
\begin{align}
  (\mVV^{m+1}\cdot\vv\nu^m,\xi^h)_{\Gm}
  -(\mV^{m+1},\xi^h)_{\Gm}&=0,
  \label{eq:nonlinear-mdr-normal}\\
  (\varkappa^{m+1},g^h)^h_{\Gm}
  +(\nabs\vvX^{m+1},\nabs(g^h\vv v^m))_{\Gm}&=0,
  \label{eq:nonlinear-mdr-curvature}\\
  (\nabs\mVV^{m+1},\nabs\vv\eta^h)_{\Gm}
  +(\mu^{m+1}\vv\nu^m,\vv\eta^h)_{\Gm}&=0.
  \label{eq:nonlinear-mdr-minimization}
\end{align}
\end{subequations}
Note that the multiplier $\mu^{m+1}$, which enforces the normal-velocity constraint, is distinct from the discrete curvature $\varkappa^{m+1}$.  

The complete gIEQ--BGN scheme is
\eqref{eq:nonlinear-ieq-common} and  \eqref{eq:nonlinear-bgn-geometry}, whereas the
complete gIEQ--MDR scheme is
\eqref{eq:nonlinear-ieq-common} and \eqref{eq:nonlinear-mdr-geometry}.  Once one
of these monolithic nonlinear systems has been solved, the fields that are needed for the next time step are pushed to the new surface by
\begin{equation}
  \lambda_{\Gamma^{m+1}}^{m+1}
  =\lambda^{m+1}\circ(\vvX^{m+1})^{-1},
  \label{eq:nonlinear-field-update}
\end{equation}
As no confusion can arise, we denoted this pushed-forward quantity
again by $\lambda^{m+1}$.

The nonlinear systems arising from \eqref{eq:nonlinear-ieq-common}, 
\eqref{eq:nonlinear-bgn-geometry} and \eqref{eq:nonlinear-ieq-common},
\eqref{eq:nonlinear-mdr-geometry}, respectively,
are solved with the help of Newton's iteration, initialized by
the solution at time $t_m$.

\subsection{Energy stability}

Using the transported field in \eqref{eq:nonlinear-field-update}, we define
\begin{equation}
  E_h^m=\frac12\norm{\lambda^m}_{L^2(\Gm)}^2,
  \qquad
  E_h^{m+1}
  =\frac12\norm{\lambda_{\Gamma^{m+1}}^{m+1}}
  _{L^2(\Gamma^{m+1})}^2.
  \label{eq:discrete-ieq-energies}
\end{equation}

We then have the following unconditional stability result, which uses ideas developed in \cite{GNZ25willmore, GNZ26splitting}.

\begin{thm}[discrete energy stability]
\label{thm:discrete-energy} Assume that $|\sigma|>0$ for all $\sigma\in\mathcal{T}^m$ and the vertex normals in \eqref{eq:vertex-normal} and \eqref{eq:implicit-vertex-normal} are well defined.  Let $(\mathscr{V}^{m+1}, \lambda^{m+1}, \mVV^{m+1}, \varkappa^{m+1})$ be a solution to the gIEQ--BGN scheme in \eqref{eq:nonlinear-ieq-common} and  \eqref{eq:nonlinear-bgn-geometry} or the gIEQ--MDR scheme in \eqref{eq:nonlinear-ieq-common} and \eqref{eq:nonlinear-mdr-geometry}, then it holds that 
\begin{equation}
  E_h^{m+1}
  +\ttau \left(\beta(\vv\nu^m)^{-1}\mV^{m+1}, \mV^{m+1}\right)_{\Gm}^h
  +\frac12\norm{\lambda^{m+1}-\lambda^m}_{L^2(\Gm)}^2
  \leq E_h^m.
  \label{eq:discrete-energy-law}
\end{equation}
\end{thm}
\begin{proof}
Set $\varphi^h=\ttau\mV^{m+1}$ in
\eqref{eq:nonlinear-ieq-v} and
$\chi^h=\ttau\lambda^{m+1}$ in
\eqref{eq:nonlinear-ieq-lambda}, and add the resulting identities.  The
$\widehat Q_\varkappa^{m+1}$ gradient terms, the $\mathcal W^m$ terms, the curvature
reaction terms and the $\widehatvvQ_{\vv p}^{m+1}$ terms cancel in
pairs.  The transport contribution vanishes because
$\mA_{\Gm}(\mVV^{m+1},\lambda^{m+1},\lambda^{m+1})=0$.  We therefore obtain
\begin{multline}
  \ttau \left(\beta(\vv\nu^m)^{-1}\mV^{m+1}, \mV^{m+1}\right)_{\Gm}^h
  +(\lambda^{m+1}-\lambda^m,\lambda^{m+1})_{\Gm}\\
  +\frac12\left(
    \nabs\vvX^{m+1}:\nabs[\vvX^{m+1}-\vv\id],
    (\lambda^{m+1})^2\right)_{\Gm}=0.
  \label{eq:nonlinear-energy-identity}
\end{multline}
The geometric estimate in \cite[(4.26)]{GNZ26splitting} yields
\begin{equation}
  \left(
    \nabs\vvX^{m+1}:\nabs[\vvX^{m+1}-\vv\id],
    (\lambda^{m+1})^2\right)_{\Gm}
  \geq
  \norm{\lambda_{\Gamma^{m+1}}^{m+1}}_{L^2(\Gamma^{m+1})}^2
  -\norm{\lambda^{m+1}}_{L^2(\Gm)}^2.
  \label{eq:nonlinear-area-inequality}
\end{equation}
For curves this follows from
$\frac12(a^2-1)\geq a-1$, while for surfaces in $\bR^3$ it follows from
the fact that the local area ratio is bounded above by
$\frac12|\nabs\vvX^{m+1}|^2$.  In both cases one also uses
$A:(A-B)\geq\frac12(|A|^2-|B|^2)$ and
$|\nabs\vv\id|^2=d-1$.  Combining
\eqref{eq:nonlinear-area-inequality} with the identity
\begin{equation*}
  (a-b)a=\frac12(a^2-b^2)+\frac12(a-b)^2,
\end{equation*}
the last two terms in \eqref{eq:nonlinear-energy-identity} are bounded below
by
\begin{equation*}
  E_h^{m+1}-E_h^m
  +\frac12\norm{\lambda^{m+1}-\lambda^m}_{L^2(\Gm)}^2.
\end{equation*}
This proves \eqref{eq:discrete-energy-law}.  The argument uses only the
common geometric--IEQ equations.  The BGN and MDR equations determine the curvature and
the tangential motion, but neither geometric formulation enters the
cancellation.
\end{proof}

\begin{rem}
To the best of our knowledge, the proposed gIEQ--MDR scheme is the first energy-stable approximation for anisotropic curvature flows with an MDR tangential velocity. The geometric--IEQ framework also accommodates the $p$-Willmore energy \cite{Aulisa2020new}
\[
E_p(\Gamma)=\frac12\int_\Gamma |\varkappa|^p\,\mathrm d\Gamma,
\qquad p\ge2\quad\mbox{and p is even},
\]
by choosing $f(\varkappa,\vec p)=|\vec p|\,\varkappa^p$ and taking $q(\varkappa,\vec p)=|\vec p|^{1/2}\varkappa^{p/2}$. The framework also applies to strongly anisotropic surface energies
with Willmore regularization,
\[
E_\varepsilon(\Gamma)
=\int_\Gamma
\left[\gamma(\vec\nu)+\frac{\varepsilon^2}{2}\varkappa^2\right]
\,\mathrm d\Gamma,
\qquad \varepsilon>0,
\]
where we take
\[
f(\varkappa,\vec p)
=2\gamma(\vec p)+\varepsilon^2|\vec p|\varkappa^2,
\qquad
q(\varkappa,\vec p)=\sqrt{f(\varkappa,\vec p)}.
\]
\end{rem}

\section{Extension to anisotropic Willmore flow}
\label{sec:aniso}

Taking the anisotropic Willmore flow as an example, we extend the
geometric--IEQ approach to $L^2$-gradient flows of energies depending on more
general geometric quantities. In particular, we let the energy depend on the anisotropic
mean curvature $\varkappa_\gamma$.  Using the evolution equation for
$\varkappa_\gamma$, we derive the corresponding geometric--IEQ formulation.  We then
present the fully discrete scheme and prove its energy stability.  The key ingredient of the extension is the evolution equation for the relevant geometric quantity.

\subsection{Geometric--IEQ formulation}
\label{subsec:aniso-ieq}

Let $\gamma:\bR^d\setminus\{\vv0\}\to\bR_{>0}$ be smooth and positively
one-homogeneous. On the surface $\Gamma(t)$ it induces the Cahn--Hoffman vector $\gamma'(\vv\nu)$, and we define
the anisotropic mean curvature by
\begin{equation}
  \varkappa_\gamma=-\nabs\cdot\gamma'(\vv\nu).
  \label{eq:anisotropic-curvature}
\end{equation}
Let $\alpha:\bR^d\setminus\{\vv0\}\to\bR_{>0}$ be smooth and positively
one-homogeneous.  We consider the modified anisotropic Willmore energy
\cite{Clarenz04wulff,Pozzi15anisoWillmore}
\begin{equation}
  W_\gamma(\Gamma)
  =\frac12\int_\Gamma \varkappa_\gamma^2  \alpha(\vv\nu)\rd\Gamma.
  \label{eq:anisotropic-energy}
\end{equation}
The density has the same homogeneity as in
\eqref{eq:general-energy}, with $\varkappa$ replaced by $\varkappa_\gamma$.  We set
\begin{equation}
  f(\varkappa_\gamma,\vv p)=\varkappa_\gamma^2\alpha(\vv p),
  \quad
  q(\varkappa_\gamma,\vv p)=\varkappa_\gamma\sqrt{\alpha(\vv p)},
  \quad
  \pqg=\sqrt{\alpha(\vv p)},
  \quad
  \npq=\frac{\varkappa_\gamma\,\alpha'(\vv p)}{2\sqrt{\alpha(\vv p)}} = \frac{q\,\alpha'(\vv p)}{2 \alpha(\vv p)},
  \label{eq:anisotropic-q}
\end{equation}
so that the auxiliary variable \eqref{eq:lambda-definition} is
$\lambda=\varkappa_\gamma\sqrt{\alpha(\vv\nu)}$, and
$W_\gamma(\Gamma)=\frac12\norm{\lambda}_{L^2(\Gamma)}^2$.

From \cite[Lemma 2.2 (2.18b)]{deckelnick2025error}, the normal evolution of the anisotropic curvature is
\begin{equation}
  \partial_t^\square\varkappa_\gamma
  =\nabs\cdot\left(\gamma''(\vv\nu)\nabs\mV\right)
   +\mV\,\lvert A\rvert_{\gamma''}^2,
  \qquad
  \lvert A\rvert_{\gamma''}^2:=\nabs\vv\nu:\gamma''(\vv\nu)\nabs\vv\nu.
  \label{eq:anisotropic-curvature-evolution}
\end{equation}
Using the fact $\partial_t^\square\varkappa_\gamma = \partial_t^\circ \varkappa_\gamma - \mVV \cdot \nabs \varkappa_\gamma$, we know the geometric-IEQ variable $\lambda$ satisfies
\begin{align}
  \partial_t^\circ \lambda
  &=\pqg\,\left(\nabs\cdot\left(\gamma''(\vv\nu)\nabs\mV\right)+\mV|A|_{\gamma''}^2\right)
    -\nabla_{\vv p}q\cdot\nabs\mV
    +\mVV\cdot\nabs q.
  \label{eq:anisotropic-q-evolution}
\end{align}
The weighted $L^2$-gradient flow of
\eqref{eq:anisotropic-energy} has normal velocity
\begin{equation}
  \mV=-\beta(\vv\nu)\left[
    \nabs\cdot\bigl(\gamma''(\vv\nu)\nabs(\lambda \pqg)\bigr)
    +\lambda \pqg\,\lvert A\rvert_{\gamma''}^2
    +\nabs\cdot(\lambda\,\npq)\right].
  \label{eq:anisotropic-normal-velocity}
\end{equation}

Testing \eqref{eq:anisotropic-q-evolution} and \eqref{eq:anisotropic-normal-velocity}, and treating transport as in
\eqref{eq:transport-skew}, gives
\begin{subequations}\label{eq:anisotropic-weak}
\begin{align}
  &\ipd{[\beta(\vv\nu)]^{-1}\mV,~\varphi}_\Gamma
  =\ipd{\gamma''(\vv\nu)\nabs(\lambda \pqg),~\nabs\varphi}_\Gamma
   -\ipd{\lambda \pqg\,\lvert A\rvert_{\gamma''}^2,~\varphi}_\Gamma\nn\\
  &\qquad
   +\frac12\ipd{\lambda^2\varkappa,~\varphi}_\Gamma
   +\ipd{\lambda\,\npq,~\nabs\varphi}_\Gamma
   \qquad\forall\varphi\in H^1(\Gt);
  \label{eq:anisotropic-weak-v}\\
  &\ipd{\partial_t^\circ\lambda,~\chi}_\Gamma
  +\frac12\ipd{\nabs\cdot\mVV,~\lambda\chi}_\Gamma
   =-\ipd{\gamma''(\vv\nu)\nabs\mV,~\nabs(\pqg\,\chi)}_\Gamma
   +\ipd{\pqg\,\chi\mV,~\lvert A\rvert_{\gamma''}^2}_\Gamma\nn\\
  &\qquad
   -\frac12\ipd{\mV,~\varkappa\lambda\chi}_\Gamma
   -\ipd{\npq\cdot\nabs\mV,~\chi}_\Gamma
   +\frac12\mA_\Gamma(\mVV,\lambda,\chi)
   \qquad\forall\chi\in H^1(\Gt).
  \label{eq:anisotropic-weak-lambda}
\end{align}
\end{subequations}
Note that the $\varkappa$ in $\frac12\ipd{\lambda^2\varkappa,~\varphi}_\Gamma$ and $\frac12\ipd{\mV,~\varkappa\lambda\chi}_\Gamma$ comes from the divergence theorem. Hence $\varkappa$ denotes the ordinary mean curvature, not the anisotropic mean curvature $\varkappa_\gamma$. The BGN and MDR formulations provide not only the tangential velocity but also the mean curvature $\varkappa$, which helps to close the system. 

With $\lambda(\cdot,0)=q(\varkappa_\gamma(\cdot,0),\vv\nu(\cdot,0))$,
choosing $(\varphi,\chi)=(\mV,\lambda)$ gives, as in
Theorem~\ref{thm:continuous-energy},
\begin{equation}
  \ddt W_\gamma(\Gt)=-\norm{[\beta(\vv \nu)]^{-\frac{1}{2}}\mV}_{L^2(\Gt)}^2 \leq 0.
  \label{eq:anisotropic-energy-law}
\end{equation}

\begin{rem}
Choosing $\alpha(\vv p)=|\vv p|$ recovers the anisotropic Willmore energy of Clarenz \cite[(6)]{Clarenz04wulff}. The choice $\alpha=\gamma$ gives the energy considered by Pozzi \cite[(2.10)]{Pozzi15anisoWillmore}; with $\beta(\vv\nu)=\gamma(\vv\nu)$, the flow \eqref{eq:anisotropic-normal-velocity} has the same normal evolution as Pozzi's anisotropic Willmore flow \cite[(3.7)--(3.8)]{Pozzi15anisoWillmore}.

For both choices, the corresponding modified anisotropic Willmore energy attains its global minimum at Wulff shapes among closed surfaces in $\bR^3$; see \cite[Theorem~4.1]{Clarenz04wulff} and \cite[Theorem~2.4]{Pozzi15anisoWillmore}, respectively.
\end{rem}

\subsection{Fully discrete geometric--IEQ and energy stability}
\label{subsec:aniso-pfem}

We use the same new-curvature, old-normal convention as in
Section~\ref{sec:pfem}.  For the signed factor
$q(\varkappa_\gamma,\vv p)=\varkappa_\gamma\sqrt{\alpha(\vv p)}$,
the curvature derivative is independent of $\varkappa_\gamma$, so that
\begin{equation}
  \widehat Q_{\varkappa_\gamma}^{m+1}= I_h^m \sqrt{\alpha(\vv v^m)},
  \qquad
  \mathcal W_{\gamma''}^m=\nabs\vv v^m: \gamma''(\vv \nu^m) \nabs\vv v^m.
  \label{eq:explicit-qkappa-coefficients2}
\end{equation}

The Hessian $\gamma''$ is also evaluated on the known surface,
whereas $\npq$ is treated implicitly using
\eqref{eq:implicit-vertex-normal}:
\begin{equation}
  \widehatvvQ_{\vv p}^{m+1}
  =\frac{\lambda^{m+1}}{2}I_h^m\left(\frac{\alpha'(\hat{\vv v}^{m+1})}{\alpha(\hat{\vv v}^{m+1})}\right).
  \label{eq:semiimplicit-qp2}
\end{equation}

The fully discrete gIEQ--BGN scheme reads as follows.  Given $\Gm$ and
$\lambda^m\in\Vh(\Gm)$, we seek
\[
  (\mV^{m+1},\lambda^{m+1},\varkappa^{m+1},
   \mVV^{m+1})
  \in[\Vh(\Gm)]^3\times[\Vh(\Gm)]^d
\]
such that, for every $(\varphi^h,\chi^h,\xi^h)\in[\Vh(\Gm)]^3$ and
$\vv\eta^h\in[\Vh(\Gm)]^d$,
\begin{subequations}\label{eq:aniso-scheme}
\begin{align}
  &([\beta(\vv\nu^m)]^{-1}\mV^{m+1},\varphi^h)_{\Gm}^h
  -(\gamma''(\vv\nu^m)[\widehat Q_{\varkappa_\gamma}^{m+1}\nabs\lambda^{m+1}
    +\lambda^{m+1}\nabs\widehat Q_{\varkappa_\gamma}^{m+1}],\nabs\varphi^h)_{\Gm}
  \nn\\
  &\quad
  +(\lambda^{m+1}\widehat Q_{\varkappa_\gamma}^{m+1},\mathcal W_{\gamma''}^m\varphi^h)_{\Gm}
  -\frac12((\lambda^{m+1})^2,\varkappa^{m+1}\varphi^h)^h_{\Gm}
  \nn\\
  &\quad
  -(\lambda^{m+1}\widehatvvQ_{\vv p}^{m+1},\nabs\varphi^h)_{\Gm}=0,
  \label{eq:aniso-scheme-v}\\
  &\left(\frac{\lambda^{m+1}-\lambda^m}{\ttau},\chi^h\right)_{\Gm}
  +\frac12\left(
    \nabs\vvX^{m+1}:\nabs\mVV^{m+1},\lambda^{m+1}\chi^h\right)_{\Gm}
  \nn\\
  &\quad
  -\frac12\mA_{\Gm}(\mVV^{m+1},\lambda^{m+1},\chi^h)
  \nn\\
  &\quad
  +(\gamma''(\vv\nu^m)\nabs\mV^{m+1},
  \widehat Q_{\varkappa_\gamma}^{m+1}\nabs\chi^h+\chi^h\nabs\widehat Q_{\varkappa_\gamma}^{m+1})_{\Gm}
  -(\widehat Q_{\varkappa_\gamma}^{m+1}\chi^h,\mathcal W_{\gamma''}^m\mV^{m+1})_{\Gm}
  \nn\\
  &\quad
  +\frac12(\mV^{m+1},\lambda^{m+1}\varkappa^{m+1}\chi^h)^h_{\Gm}
  +(\widehatvvQ_{\vv p}^{m+1},\nabs\mV^{m+1}\chi^h)_{\Gm}=0,
  \label{eq:aniso-scheme-lambda}\\
  &(\mVV^{m+1}\cdot\vv \nu^m,\xi^h)^h_{\Gm}
  -(\mV^{m+1},\xi^h)^h_{\Gm}=0,
  \label{eq:aniso-scheme-normal}\\
  &(\varkappa^{m+1}\vv\nu^m,\vv\eta^h)^h_{\Gm}
  +(\nabs\vvX^{m+1},\nabs\vv\eta^h)_{\Gm}=0.
  \label{eq:aniso-scheme-curvature}
\end{align}
\end{subequations}
The mean curvature $\varkappa^{m+1}$ is determined by
the BGN curvature equation \eqref{eq:aniso-scheme-curvature}. For the gIEQ--MDR variant, \eqref{eq:aniso-scheme-normal} and \eqref{eq:aniso-scheme-curvature} are replaced by \eqref{eq:nonlinear-mdr-geometry}.

Under the assumptions of Theorem~\ref{thm:discrete-energy}, testing
\eqref{eq:aniso-scheme-v} and \eqref{eq:aniso-scheme-lambda} with
$\ttau\mV^{m+1}$ and $\ttau\lambda^{m+1}$ gives the same cancellation.
The symmetry of $\gamma''(\vv\nu^m)$ pairs the anisotropic gradient
terms, and the remaining terms cancel as before.  Thus
\eqref{eq:discrete-energy-law} holds for both BGN and MDR variants, with
$E_h^m=\frac12\norm{\lambda^m}_{L^2(\Gm)}^2$.

\section{Geometric--IEQ formulation for the \texorpdfstring{$H^{-1}$}{H-1}-gradient flow}
\label{sec:h-1}
\label{subsec:h-1}

So far we have considered geometric--IEQ for the $L^2$-gradient flows. We now extend the approach to $H^{-1}$-gradient flows, which arise in many applications. A typical example is surface diffusion. Compared with $L^2$-gradient flows, $H^{-1}$-gradient flows not only dissipate the energy but also conserve the enclosed mass. We construct geometric-IEQ PFEMs for $H^{-1}$-gradient flows that preserve both geometric properties.

\subsection{Weak formulation and conservation law}

For simplicity, we consider the energy \eqref{eq:general-energy} of Section~\ref{sec:mf}. The equation for the auxiliary variable is unchanged; only the normal velocity is modified. We therefore present only the differences, referring to Sections~\ref{sec:mf} and~\ref{sec:pfem} for the common parts. Introducing $w$, with the sign convention of \eqref{eq:ieq-normal-velocity}, the weighted $H^{-1}$-gradient flow reads 
\begin{equation}
  \mV=-\nabs\cdot\bigl(\beta(\vv\nu)\,\nabs w\bigr),
  \qquad
  w=-\Delta_s(\lambda \pq)-\lambda \pq\lvert\nabs\vv\nu\rvert^2 -\nabs\cdot(\lambda\,\npq) \quad\mbox{on}\quad\Gt.
  \label{eq:h-1-flow}
\end{equation}
When the density in \eqref{eq:general-energy} is simply
$f(\varkappa,\vv\nu) = 2|\vv\nu|=2$ on $\Gamma$, then \eqref{eq:h-1-flow} 
with $\beta(\vv\nu)=1$ reduces to surface diffusion, with $w = \varkappa$,
which is the well-known $H^{-1}$-gradient flow of surface area.

We replace \eqref{eq:general-weak-v} by: find
$(\mV,w)\in H^1(\Gt)\times H^1(\Gt)$ such that
\begin{subequations}\label{eq:h-1-weak}
\begin{align}
  &\ipd{\mV,~\varpi}_\Gamma
  -\ipd{\beta(\vv\nu)\,\nabs w,~\nabs\varpi}_\Gamma=0
  \qquad\forall\varpi\in H^1(\Gt);
  \label{eq:h-1-weak-v}\\
  &\ipd{w,~\varphi}_\Gamma
  =\ipd{\nabs(\lambda \pq),~\nabs\varphi}_\Gamma
   -\ipd{\lambda\pq\lvert\nabs\vv\nu\rvert^2,~\varphi}_\Gamma\nn\\
  &\qquad\qquad\qquad
   +\frac12\ipd{\lambda^2\varkappa,~\varphi}_\Gamma
   +\ipd{\lambda\npq,~\nabs\varphi}_\Gamma
  \qquad\forall\varphi\in H^1(\Gt).
  \label{eq:h-1-weak-w}
\end{align}
\end{subequations}
These equations are coupled with \eqref{eq:general-weak-lambda} and
either \eqref{eqn:general-weakBGN} or \eqref{eqn:general-weakMDR}.
Testing with $(\varpi,\varphi,\chi)=(w,\mV,\lambda)$ gives
\begin{equation}
  \ddt\left(\frac12\norm{\lambda}_{L^2(\Gt)}^2\right)
  =-\norm{[\beta(\vv\nu)]^{\frac12}\nabs w}_{L^2(\Gt)}^2\leq0.
  \label{eq:h-1-energy}
\end{equation}
In addition, choosing $\varpi=1$ yields
\begin{equation}
  \ddt{\rm vol}(\Gamma(t))=\int_{\Gt}\mV\rd\Gamma=0,
  \label{eq:h-1-mass}
\end{equation}
where ${\rm vol}(\Gamma(t))$ represents the volume enclosed by $\Gt$.

\subsection{Structure-preserving geometric-IEQ}

In addition to the unknowns in Section~\ref{subsec:nonlinear-pfem}, we
seek $w^{m+1}\in\Vh(\Gm)$.  The discrete geometric--IEQ equations consist of
\eqref{eq:nonlinear-ieq-lambda} and, for every
$(\varpi^h,\varphi^h)\in\Vh(\Gm)^2$,
\begin{subequations}\label{eq:h-1-ieq-common}
\begin{align}
  &(\mV^{m+1},\varpi^h)^h_{\Gm}
  -(\beta(\vv\nu^m)\nabs w^{m+1},\nabs\varpi^h)_{\Gm}=0,
  \label{eq:h-1-fd-v}\\
  &(w^{m+1},\varphi^h)^h_{\Gm}
  -(\widehat Q_\varkappa^{m+1}\nabs\lambda^{m+1}
    +\lambda^{m+1}\nabs\widehat Q_\varkappa^{m+1},\nabs\varphi^h)_{\Gm}
  \nn\\
  &\quad +(\lambda^{m+1}\widehat Q_\varkappa^{m+1},
    \mathcal W^m\varphi^h)_{\Gm}
  -\frac12((\lambda^{m+1})^2,
    \varkappa^{m+1}\varphi^h)^h_{\Gm}
  \nn\\
  &\quad -(\lambda^{m+1}\widehatvvQ_{\vv p}^{m+1},
    \nabs\varphi^h)_{\Gm}=0.
  \label{eq:h-1-fd-w}
\end{align}
\end{subequations}
Thus the only change to \eqref{eq:nonlinear-ieq-v} is the replacement
of $[\beta(\vv\nu^m)]^{-1}\mV^{m+1}$ by $w^{m+1}$, together with
the additional equation \eqref{eq:h-1-fd-v}.

To preserve the enclosed mass, we use the time-averaged normal $\vv\nu^{m+\frac12}$ of
\cite{BZ21SPFEM}. For gIEQ--BGN we retain \eqref{eq:nonlinear-bgn-curvature} and replace \eqref{eq:nonlinear-bgn-normal} by
\begin{equation}
  (\mVV^{m+1}\cdot\vv\nu^{m+\frac12},\xi^h)^h_{\Gm}
  -(\mV^{m+1},\xi^h)^h_{\Gm}=0
  \qquad\forall\xi^h\in\Vh(\Gm),
  \label{eq:sp-bgn-normal}
\end{equation}
where $\vec\nu^{m+\frac12}$ satisfies 
\begin{equation}
  {\rm vol}(\Gamma^{m+1})-{\rm vol}(\Gamma^m)
  =\ttau(\mVV^{m+1}\cdot\vv\nu^{m+\frac12},1)^h_{\Gm}.\label{eq:volumep}
\end{equation}
For gIEQ--MDR we retain \eqref{eq:nonlinear-mdr-curvature} and
\eqref{eq:nonlinear-mdr-minimization}, and replace
\eqref{eq:nonlinear-mdr-normal} by \eqref{eq:sp-bgn-normal}. The surface updates are as in \eqref{eq:nonlinear-surface-map} and \eqref{eq:nonlinear-field-update}.

\begin{thm}[energy stability and mass conservation]
\label{thm:h-1-discrete-energy}
 Assume that $|\sigma|>0$ for all $\sigma\in\mathcal{T}^m$ and the vertex normals in \eqref{eq:vertex-normal} and \eqref{eq:implicit-vertex-normal} are well defined.  Let $(\mathscr{V}^{m+1}, \lambda^{m+1}, w^{m+1}, \vec X^{m+1}, \varkappa^{m+1})$ be a solution to the gIEQ--BGN scheme in \eqref{eq:h-1-ieq-common}, \eqref{eq:nonlinear-ieq-lambda}, \eqref{eq:sp-bgn-normal} and \eqref{eq:nonlinear-bgn-curvature} or the gIEQ--MDR scheme in \eqref{eq:h-1-ieq-common}, \eqref{eq:nonlinear-ieq-lambda}, \eqref{eq:sp-bgn-normal}, \eqref{eq:nonlinear-mdr-curvature} and \eqref{eq:nonlinear-mdr-minimization},  then the discrete solutions satisfy
\begin{equation}
  E_h^{m+1}
  +\ttau\norm{[\beta(\vv\nu^m)]^{\frac12}\nabs w^{m+1}}_{L^2(\Gm)}^2
  +\frac12\norm{\lambda^{m+1}-\lambda^m}_{L^2(\Gm)}^2
  \leq E_h^m,
  \label{eq:h-1-discrete-energy}
\end{equation} and
\begin{equation}
  {\rm vol}(\Gamma^{m+1})={\rm vol}(\Gamma^m)
  \label{eq:h-1-discrete-mass}
\end{equation}
with $E_h^m$ defined in \eqref{eq:discrete-ieq-energies}, and  ${\rm vol}(\Gm)$ is the volume enclosed by the discrete surface $\Gm$.
\end{thm}
\begin{proof}
Choose $\varphi^h=\ttau\mV^{m+1}$ in \eqref{eq:h-1-fd-w} and
$\chi^h=\ttau\lambda^{m+1}$ in \eqref{eq:nonlinear-ieq-lambda}.
The cancellation in Theorem~\ref{thm:discrete-energy} is unchanged,
with the dissipation term replaced by
$\ttau(w^{m+1},\mV^{m+1})^h_{\Gm}$.  Testing
\eqref{eq:h-1-fd-v} with $\varpi^h=w^{m+1}$ identifies this term as
$\ttau(\beta(\vv\nu^m)\nabs w^{m+1},\nabs w^{m+1})_{\Gm}$.
The estimate \eqref{eq:nonlinear-area-inequality} then gives
\eqref{eq:h-1-discrete-energy}.

Choosing $\xi^h=1$ in \eqref{eq:sp-bgn-normal}, $\varpi^h=1$ in \eqref{eq:h-1-fd-v}, and recalling \eqref{eq:volumep} gives the exact mass preservation. 
\end{proof}

\section{Numerical results}
\label{sec:num-general}

In this section, we present numerical experiments to illustrate the performance of the proposed gIEQ--BGN and gIEQ--MDR methods.  These include a convergence test, investigations of the numerical energy stability and mass conservation, and morphological evolutions. 

We use a uniform time step size $\tau$. For polygonal curves and polyhedral surfaces, we denote by $h$ the maximal edge length of $\Gamma^0$. We implement the gIEQ--BGN and gIEQ--MDR schemes with piecewise linear parametric finite elements via NGSolve \cite{Schoeberl14ngsolve}. The nonlinear systems are solved by Newton's method with tolerance $10^{-10}$.  We use the gIEQ--BGN schemes for curves, since the BGN formulation leads to an asymptotic equidistribution of the mesh points \cite{BGN07}, and the gIEQ--MDR schemes for surfaces unless otherwise stated.

For a smooth initial surface $\Gamma(0)$ with a given anisotropic energy, the initial geometric quantities $\varkappa$, $\vv\nu$, and $\varkappa_\gamma$ can be computed analytically, and the initial auxiliary variable is the nodal interpolation $\lambda^0=I_h^0\,q(\varkappa(\cdot,0),\vv\nu(\cdot,0))$. Here $\varkappa$ is replaced by $\varkappa_\gamma$ for the anisotropic Willmore energies. 

Most of the examples employ the anisotropies of BGN type
\cite{BGN08Ani,BGN08ani3d}
\begin{equation}
  \gamma(\vv p)=\sum_{\ell=1}^{L}\sqrt{\vv p^TG_\ell\,\vv p},
  \label{eq:num-gamma-G}
\end{equation}
where $G_\ell\in\bR^{d\times d}$, $\ell=1,\ldots,L$, are symmetric and positive
definite.  In the case $L=1$, the Wulff shape
$\mathcal W=\{\vv x\in\bR^d:\vv x^TG_1^{-1}\vv x\leq1\}$ is an ellipse
($d=2$) or an ellipsoid ($d=3$). More specifically, for $\Gamma(t)=R\,\partial\mathcal W$, the anisotropic mean curvature is $\varkappa_\gamma=-(d-1)/R$.

\subsection{Convergence tests}
\label{subsec:num-convergence}

We first test the convergence orders in space and in time with the help of exact self-similar solutions. We take $L=1$ in \eqref{eq:num-gamma-G}, with
$G_1=\operatorname{diag}(1,4)$ for $d=2$ and $G_1=\operatorname{diag}(1,1,4)$ for
$d=3$, so that the Wulff shape $\mathcal W$ is a $1 \times 2$ ellipse or a $1\times 1\times 2$ ellipsoid.  We consider the
following two flows:
\begin{itemize}
\item the anisotropic mean curvature flow (MCF), i.e.\ the $L^2$-gradient
  flow \eqref{eq:ieq-normal-velocity} with $f=\gamma$,
  $q=\sqrt{\gamma}$ and $\beta=\gamma$.  Its exact solution is
  $\Gamma(t)=R(t)\,\partial\mathcal W$ with $R(t)=\sqrt{1-(d-1)t}$;

\item the anisotropic Willmore flow of Pozzi type \cite[Example~3.1]{Pozzi15anisoWillmore} for curves, i.e.\
  \eqref{eq:anisotropic-normal-velocity} with $\alpha=\beta=\gamma$.  Its
  exact solution is $\Gamma(t)=R(t)\,\partial\mathcal W$ with
  $R(t)=(1+2t)^{1/4}$.
\end{itemize}
For anisotropic MCF, we consider the energy $\frac{1}{2}\int_{\Gamma}\gamma(\vec\nu)\rd\Gamma$. Consequently, the resulting gradient flow evolves at half the speed of the standard anisotropic MCF.

For curves, the initial polygon is $\vvX^0(\rho_j)=(\cos\theta_j,
2\sin\theta_j)^T$, $\theta_j=2\pi j/N$, with $N=2^k$, $k=4,\ldots,10$.
For surfaces, the initial meshes are obtained by
refining an icosahedron $r$ times at the edge midpoints, projecting the
vertices onto the unit sphere, and mapping them by
$\vv x\mapsto\operatorname{diag}(1,1,2)\,\vv x$ with $r=1,\ldots,5$.  The time step sizes are $\tau=2^{-k}$ with
$k=6,\ldots,14$ for curves and $k=4,\ldots,14$ for surfaces.

Let $M=T/\tau$ with $T=0.125$, let $\pi$ be the closest-point projection onto
$\Gamma(T)$, and let $\varkappa_\gamma(\cdot,T),\vv\nu(\cdot,T)$ be the exact geometric quantities on $\Gamma(T)$. We define the errors on $\Gamma^M$ by
\begin{equation}
  e_x(\vv x)=|\vv x-\pi(\vv x)|,
  \qquad
  e_\lambda(\vv x)=\lambda^M(\vv x)-q(\varkappa_\gamma(\pi(\vv x),T),\vv\nu(\pi(\vv x),T)),
  \qquad \vv x\in\Gamma^M,
  \label{eq:num-errors}
\end{equation}
and consider the $L^2(\Gamma^M)$-norms, denoted by $\norm{\cdot}_{L^2(\Gamma(T))}$ for brevity. 

\begin{figure}[htp!]
\centering
\includegraphics[width=0.4\textwidth]{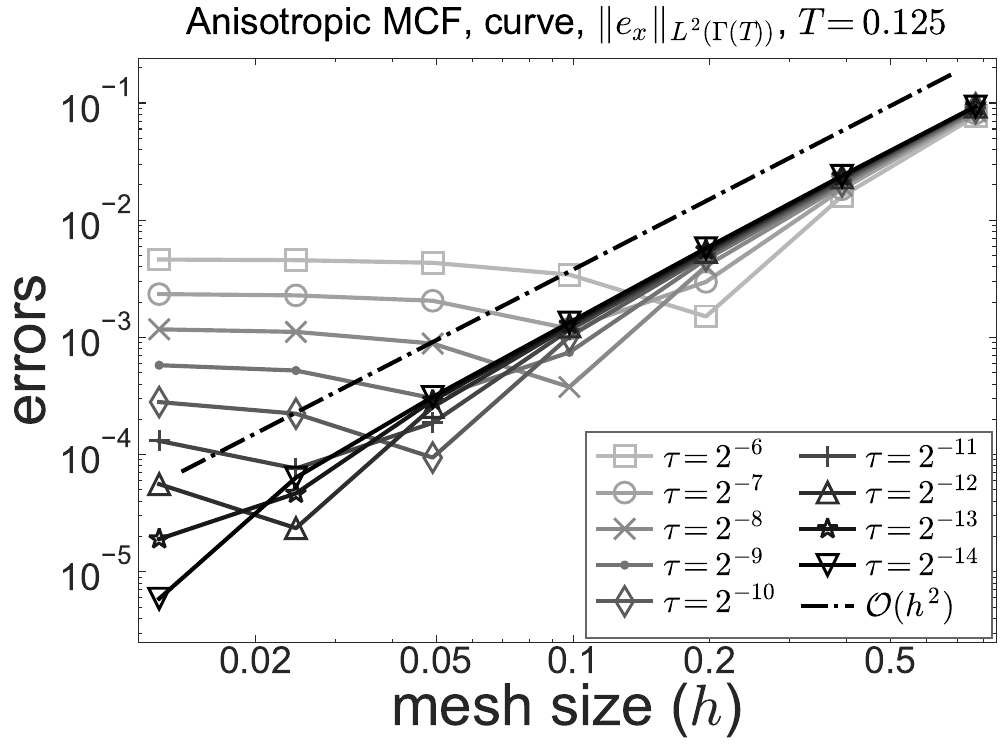}\qquad
\includegraphics[width=0.4\textwidth]{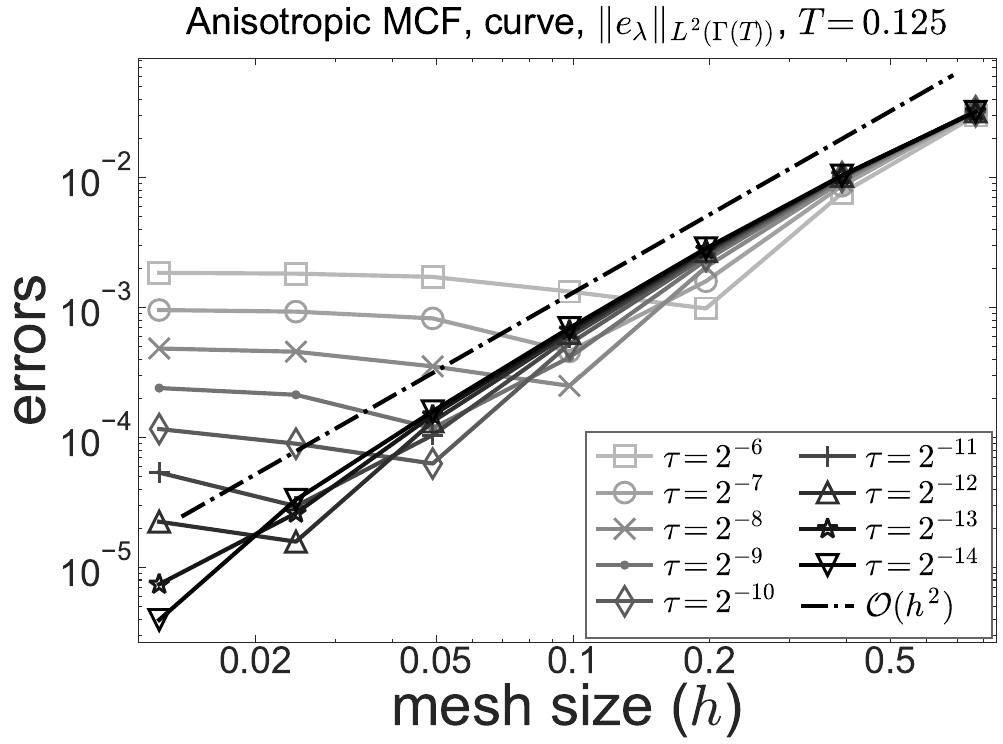}\\[0.4em]
\includegraphics[width=0.4\textwidth] {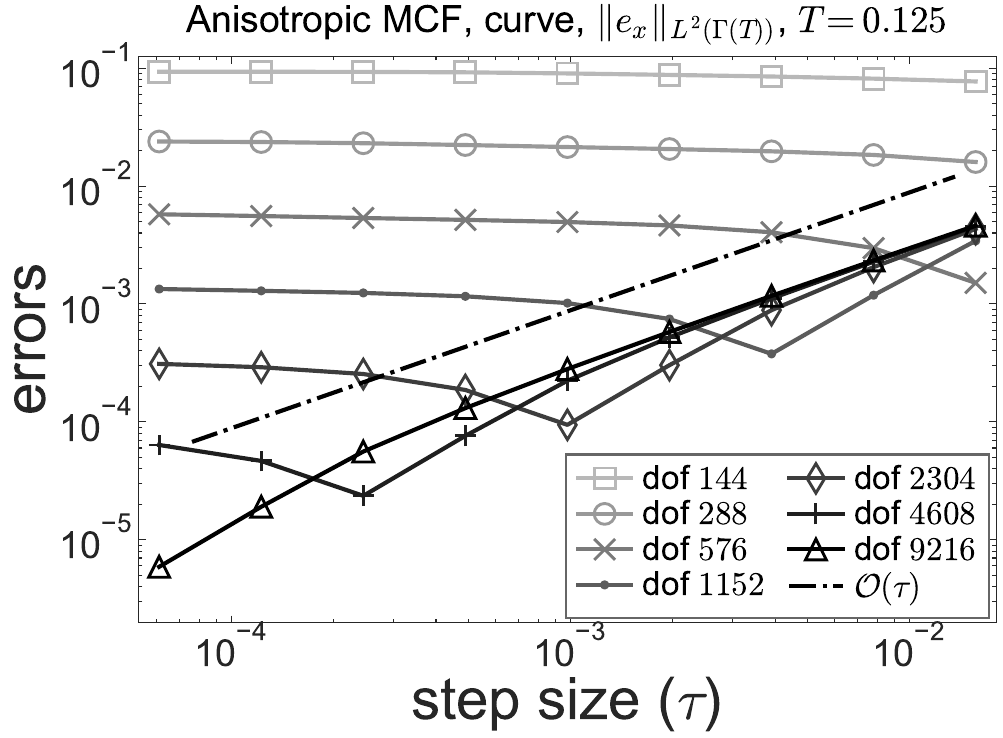}\qquad
\includegraphics[width=0.4\textwidth]{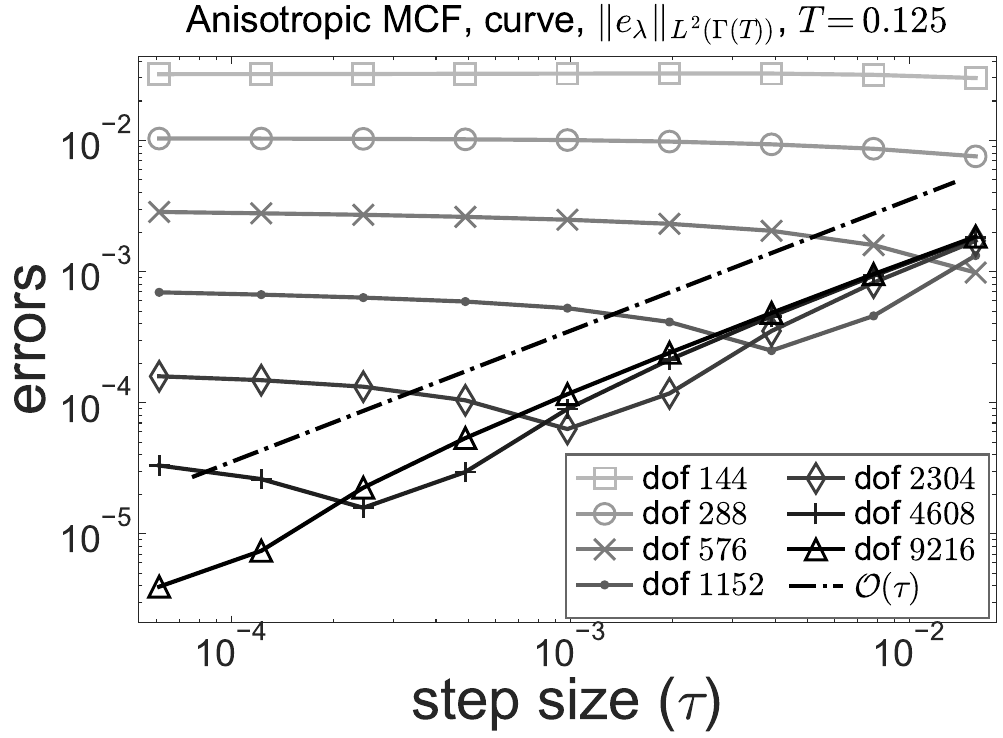}\\[0.4em]
\includegraphics[width=0.4\textwidth]{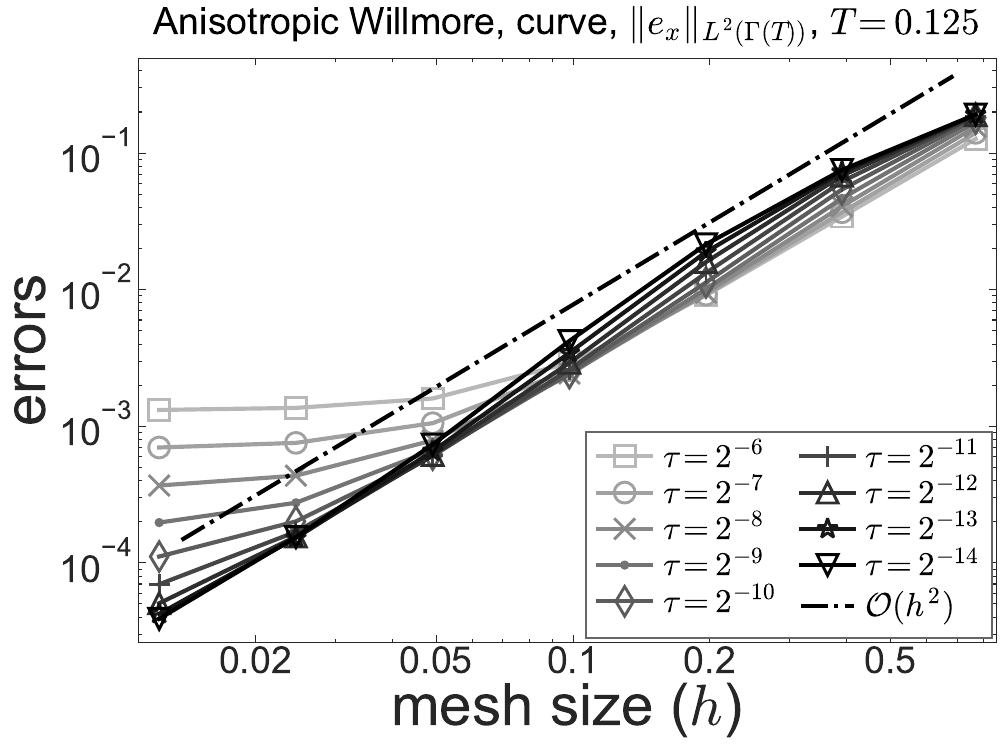}\qquad
\includegraphics[width=0.4\textwidth]{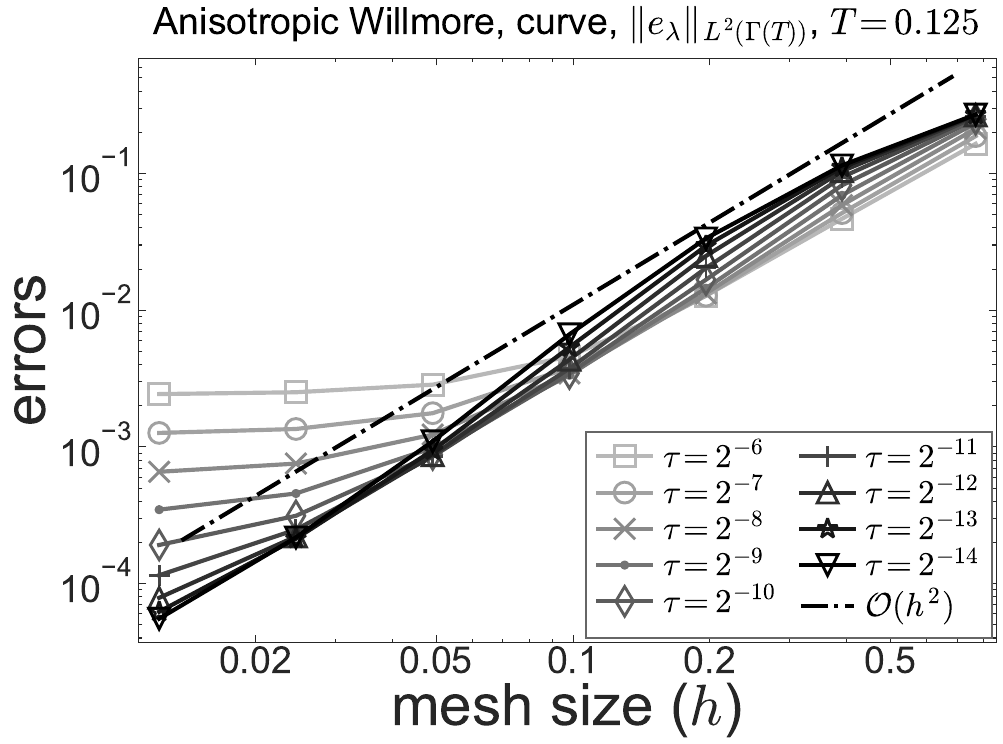}\\[0.4em]
\includegraphics[width=0.4\textwidth]{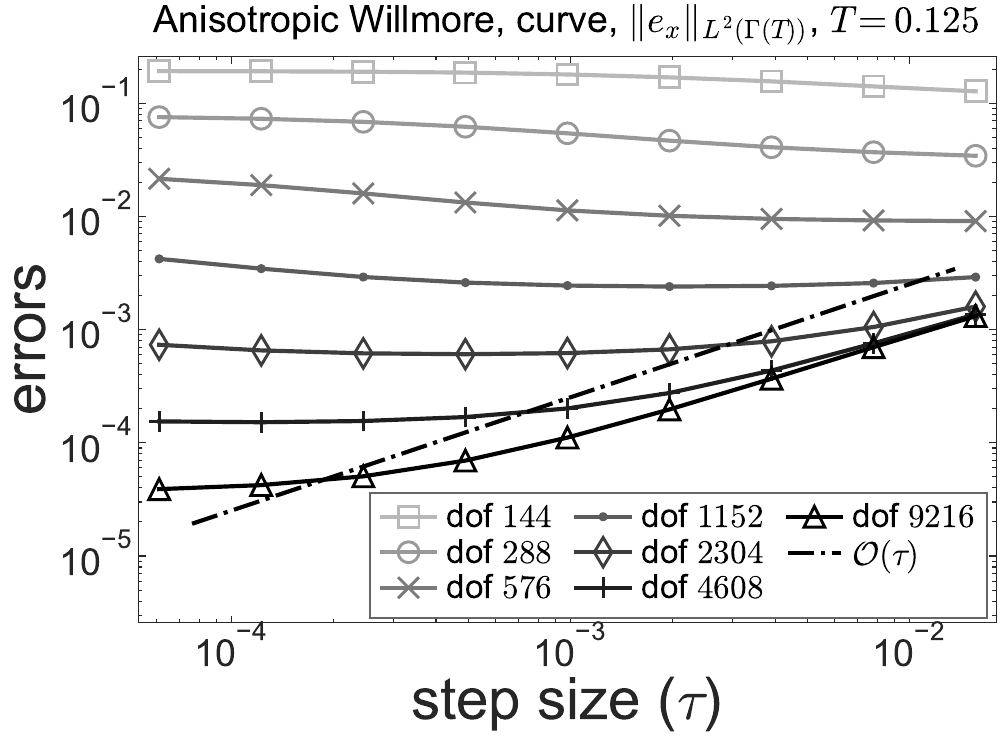}\qquad
\includegraphics[width=0.4\textwidth]{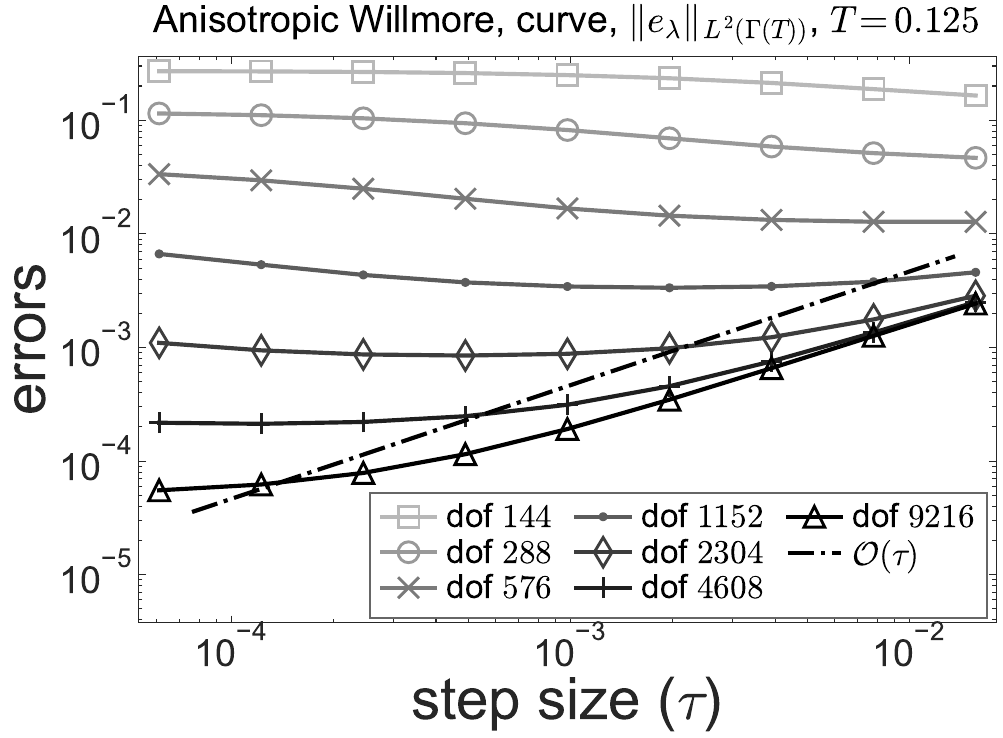}
\caption{Numerical errors of the gIEQ--BGN schemes for curves at $T=0.125$.
The first four correspond to the anisotropic MCF, and the last four
to the anisotropic Willmore flow of Pozzi.  In each case, the upper row shows
$\norm{e_x}$ (left) and $\norm{e_\lambda}$ (right) versus $h$ for fixed
$\tau$, and the lower row shows them versus $\tau$ for fixed meshes.}
\label{fig:num-conv-curves}
\end{figure}

In Fig.~\ref{fig:num-conv-curves} we report the errors at $T=0.125$ for a
sequence of meshes and a sequence of time steps. In the plots versus $h$, we can observe two regions: a region where the spatial discretization error dominates, matching the $\mathcal O(h^2)$ order of convergence, and a region, with small mesh size, where the temporal discretization error dominates.  The same situation applies to the plots versus $\tau$, with the $\mathcal O(\tau)$ order of convergence. 

\begin{figure}[htp!]
\centering
\includegraphics[width=0.4\textwidth]{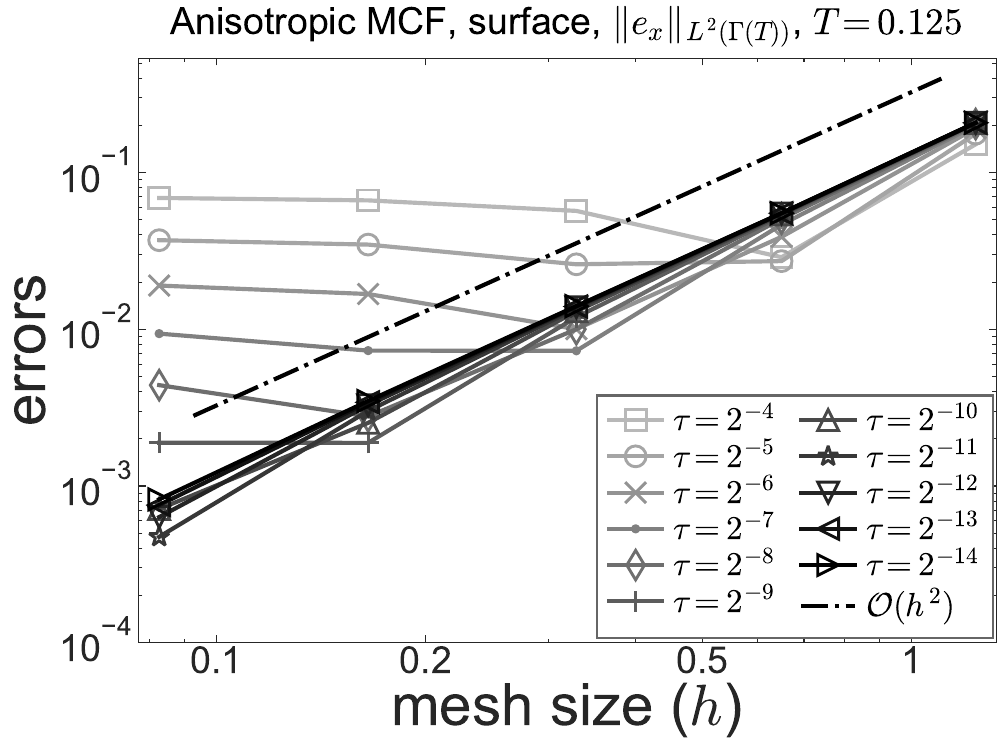}\qquad
\includegraphics[width=0.4\textwidth]{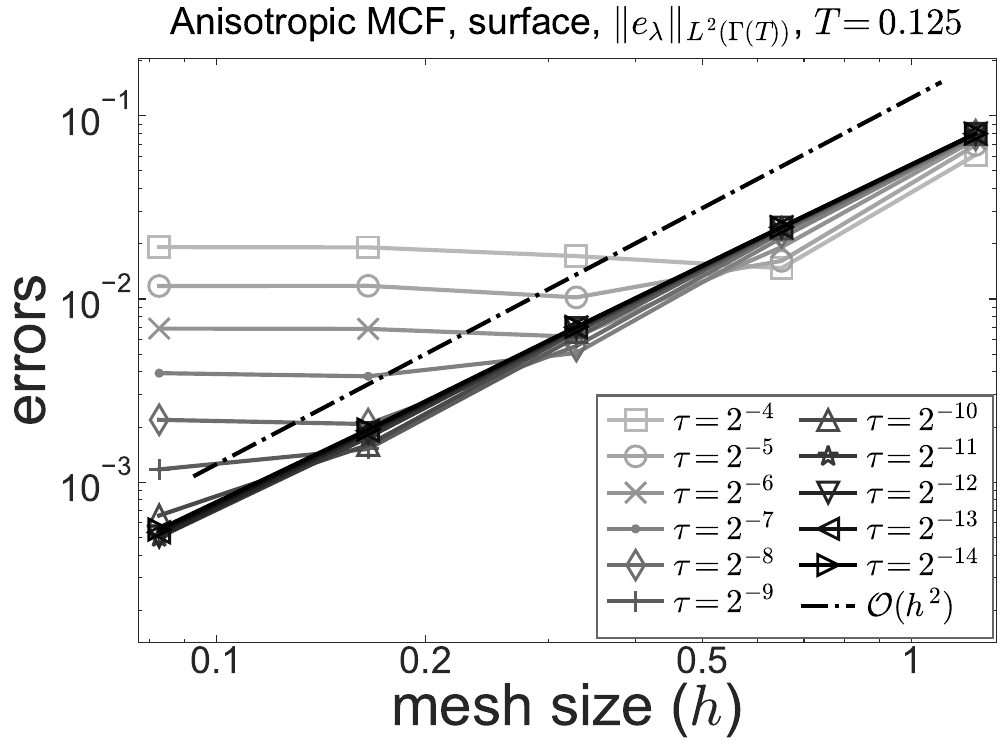}\\[0.4em]
\includegraphics[width=0.4\textwidth]{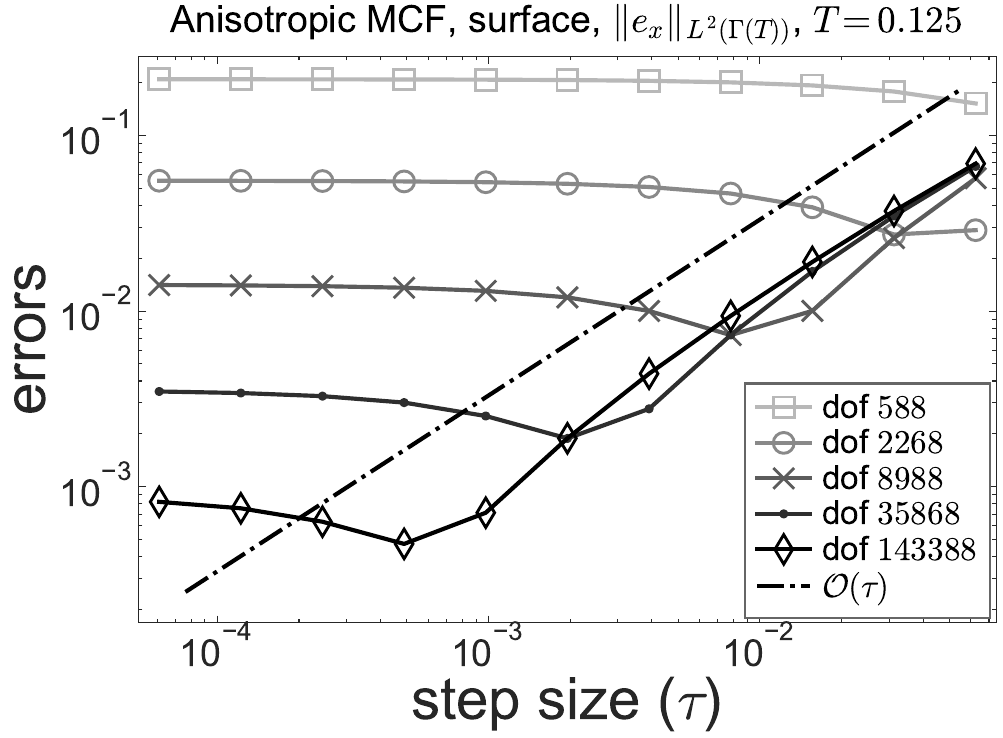}\qquad
\includegraphics[width=0.4\textwidth]{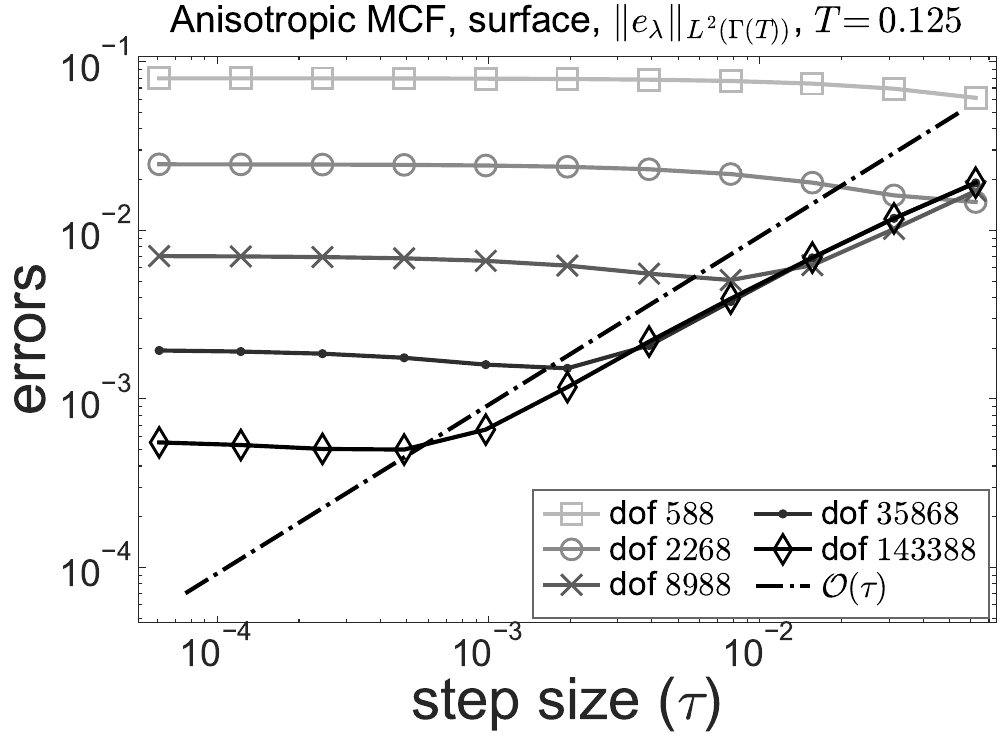}
\caption{Numerical errors of the gIEQ--MDR scheme for the anisotropic MCF of
surfaces at $T=0.125$.  The upper row shows $\norm{e_x}$ (left) and
$\norm{e_\lambda}$ (right) versus $h$ for fixed $\tau$, and the lower row shows
them versus $\tau$ for fixed meshes.}
\label{fig:num-conv-surface-mcf}
\end{figure}

In Fig.~\ref{fig:num-conv-surface-mcf} we report the corresponding errors for surfaces, and the same situation applies. We observe the convergence orders $\mathcal O(h^2)$ in space and $\mathcal O(\tau)$ in time for both errors in all scenarios.  Therefore, in the following examples we choose $\tau=\mathcal O(h^2)$ to balance the spatial and temporal discretization errors.

\subsection{Energy stability and mass conservation}
\label{subsec:num-indicators}

Next we test the energy stability and the mass conservation. In addition, we present the number of Newton iterations and the IEQ defect between the IEQ variable $\lambda$ and the geometric energy variable  $q$.  To this end, we introduce the following indicators:
the normalized IEQ energy $E_h^m/E_h^0$ with $E_h^m$ given in
\eqref{eq:discrete-ieq-energies}; the normalized geometric energy
$E(\Gamma^m)/E(\Gamma^0)$ defined on $\Gamma^m$; the
enclosed area $A^m$ for curves, or the enclosed volume $V^m$ for surfaces; the number of Newton iterations in each time step; and the IEQ defect
\begin{equation}
  \norm{\lambda^m-q_h^m}_{L^2(\Gamma^m)},
  \qquad
  q_h^m=I_h^m\,q(\varkappa_h^m,\vv v^m),
  \label{eq:num-defect}
\end{equation}
between the auxiliary variable $\lambda$ and its geometric counterpart $q$.  Here
$\vv v^m$ is the vertex normal \eqref{eq:vertex-normal}. For the anisotropic Willmore energy, $\varkappa_h$ is replaced by the anisotropic mean curvature
$\varkappa_{\gamma,h}^m\in\Vh(\Gm)$ given by
\begin{equation}
  \ipd{\varkappa_{\gamma,h}^m,~\xi^h}^h_{\Gm}
  =-\ipd{\gamma''(\vv\nu^m):\nabs\vv v^m,~\xi^h}_{\Gm}
  \qquad\forall\xi^h\in\Vh(\Gm).
  \label{eq:num-kappa-gamma}
\end{equation}
We present both the curves and the surfaces.

\textbf{Curves.}  The initial curve is the flower
$\vvX^0(\rho_j)=r(\theta_j)(\cos\theta_j,\sin\theta_j)^T$ with
$r(\theta)=1+0.2\cos(5\theta)$ and $\theta_j=2\pi j/N$, where $N=128$ and
$N=256$.  We take $L=1$ and $G_1=\operatorname{diag}(1,4)$ in \eqref{eq:num-gamma-G}, and
$\beta=\gamma$, and consider the anisotropic MCF and the anisotropic surface
diffusion (SD) of
$E(\Gamma)=\frac12\int_\Gamma\gamma(\vv\nu)\rd\Gamma$ by the gIEQ--BGN schemes
of Sections~\ref{sec:pfem} and \ref{sec:h-1}.  The time step size is
$\tau=h^2/10$. 

\textbf{Surfaces.} The initial surface is the unit sphere, discretized by the icosahedral meshes of Section~\ref{subsec:num-convergence} with $N=2562$ and $N=10242$. We consider the $H^{-1}$-gradient flow of the anisotropic Willmore energy \eqref{eq:anisotropic-energy} of Clarenz, with $\alpha(\vv p)=|\vv p|$, $\beta=1$, $L=1$, and $G_1=\operatorname{diag}(1,\frac14,\frac14)$ in \eqref{eq:num-gamma-G}. We employ the gIEQ--BGN scheme of Section~\ref{sec:h-1} with $\tau=h^2/2$, yielding $13N=33306$ and $133146$ degrees of freedom, respectively. 

\begin{figure}[htp!]
\centering
\includegraphics[width=0.44\textwidth]{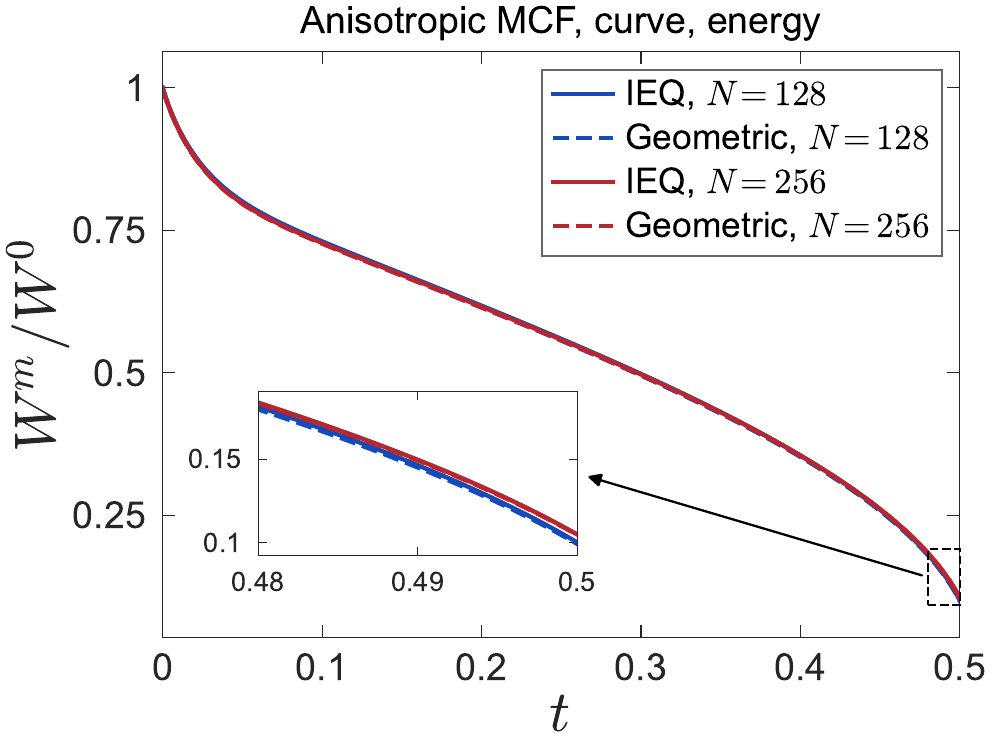}\qquad
\includegraphics[width=0.44\textwidth]{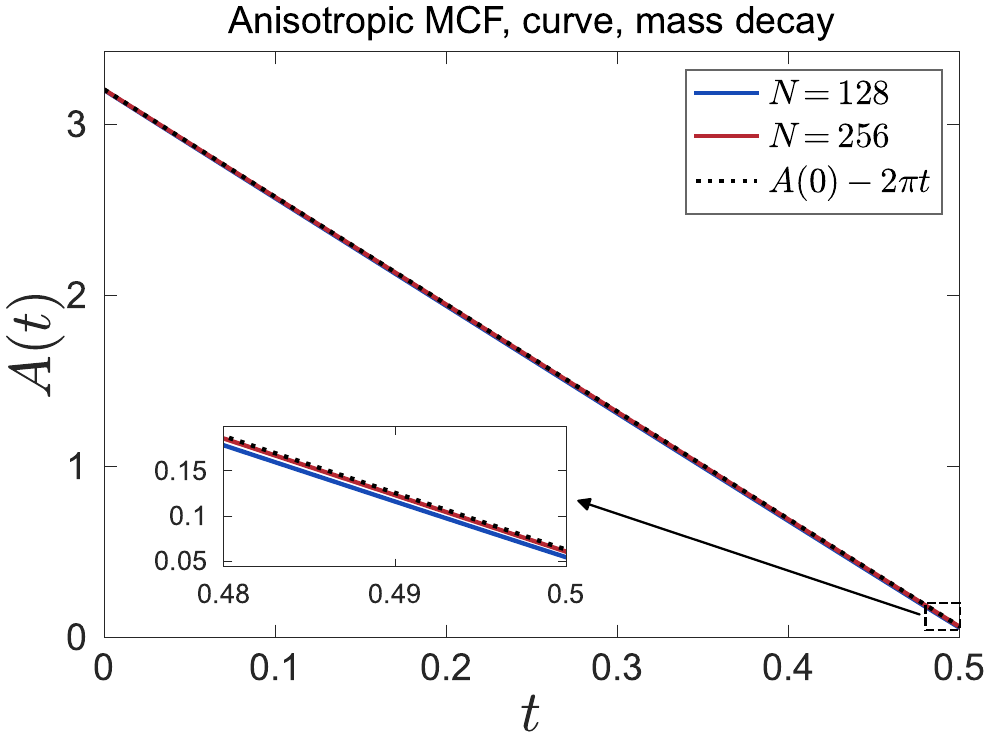}
\includegraphics[width=0.44\textwidth]{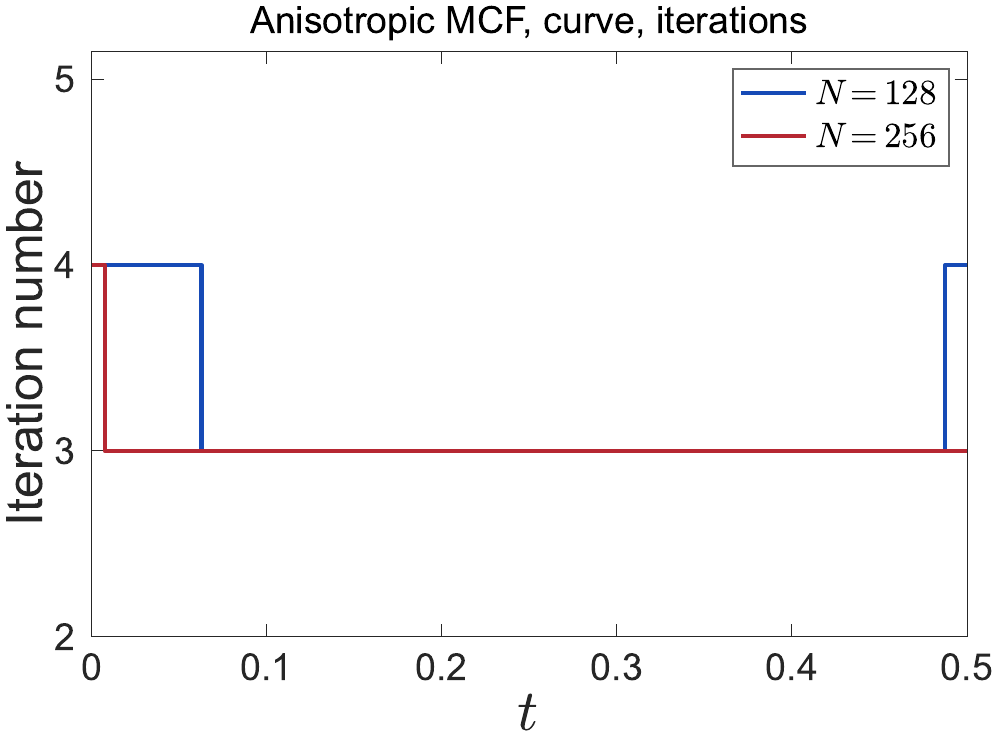}\qquad
\includegraphics[width=0.44\textwidth]{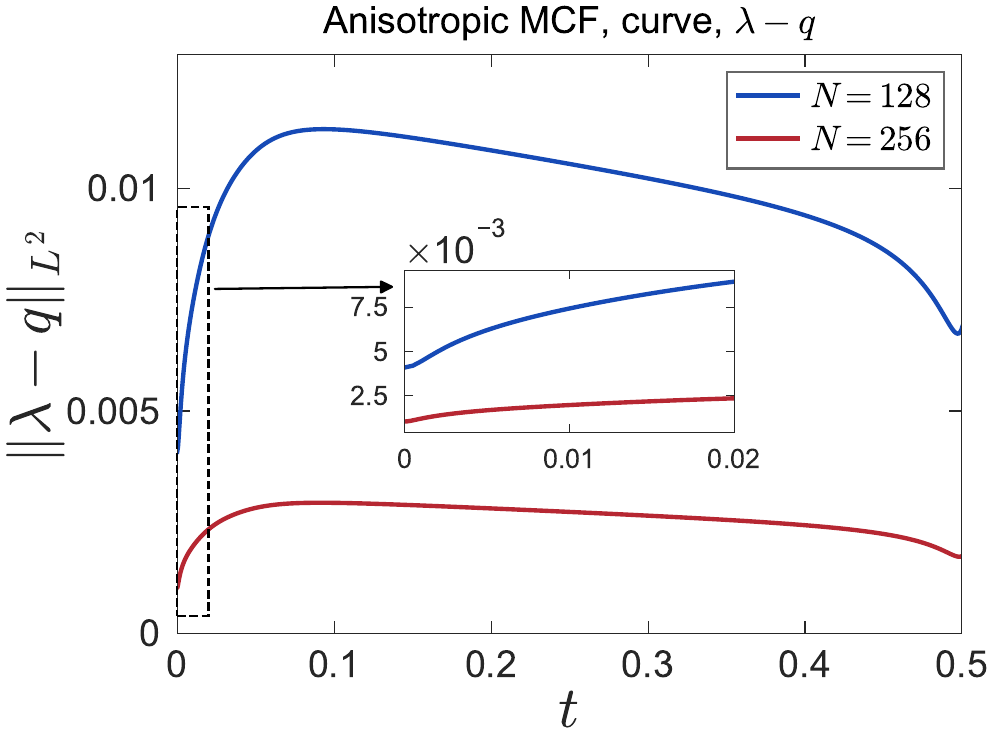}
\caption{Anisotropic MCF of the flower curve.  Upper row: the normalized IEQ and geometric energies (left) and the enclosed area $A^m$, where the dotted line is $A(0)-2\pi t$ (right). Lower row: the number of Newton iterations (left) and the IEQ defect
\eqref{eq:num-defect} (right).}
\label{fig:num-ind-l2}
\end{figure}

\begin{figure}[htp!]
\centering
\includegraphics[width=0.4\textwidth]{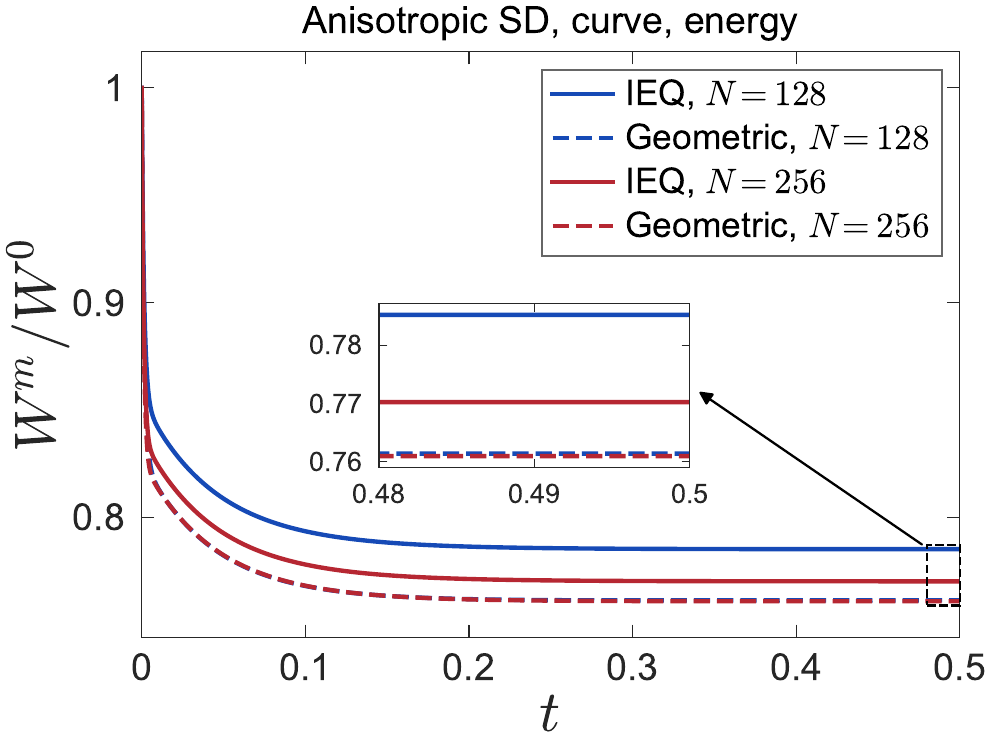}\qquad
\includegraphics[width=0.4\textwidth]{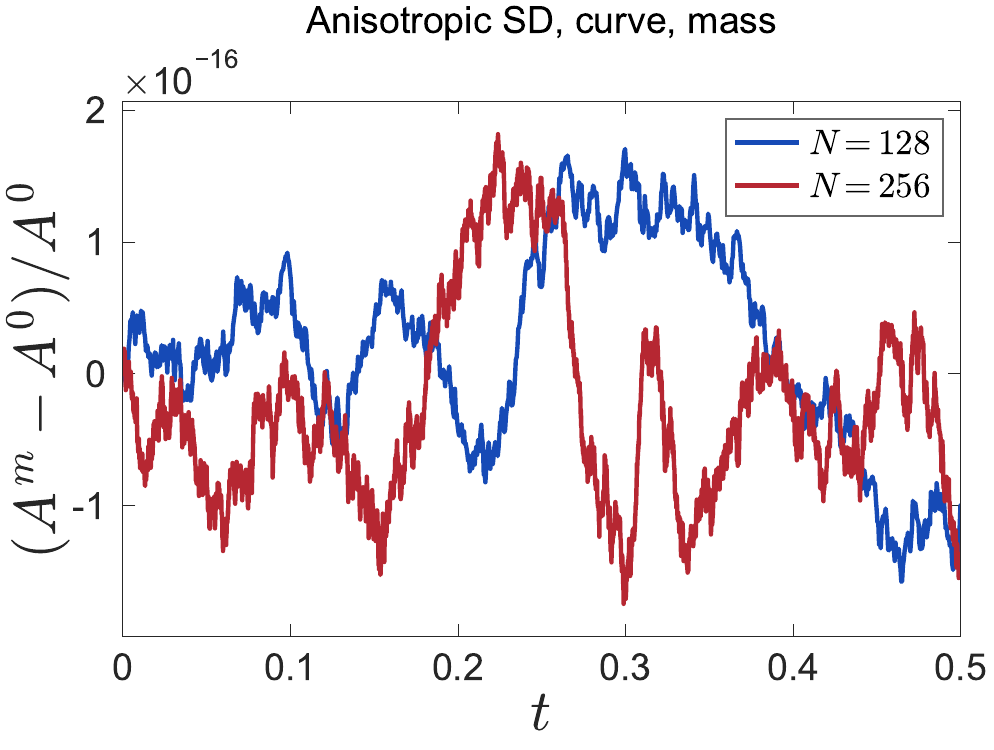}
\includegraphics[width=0.4\textwidth]{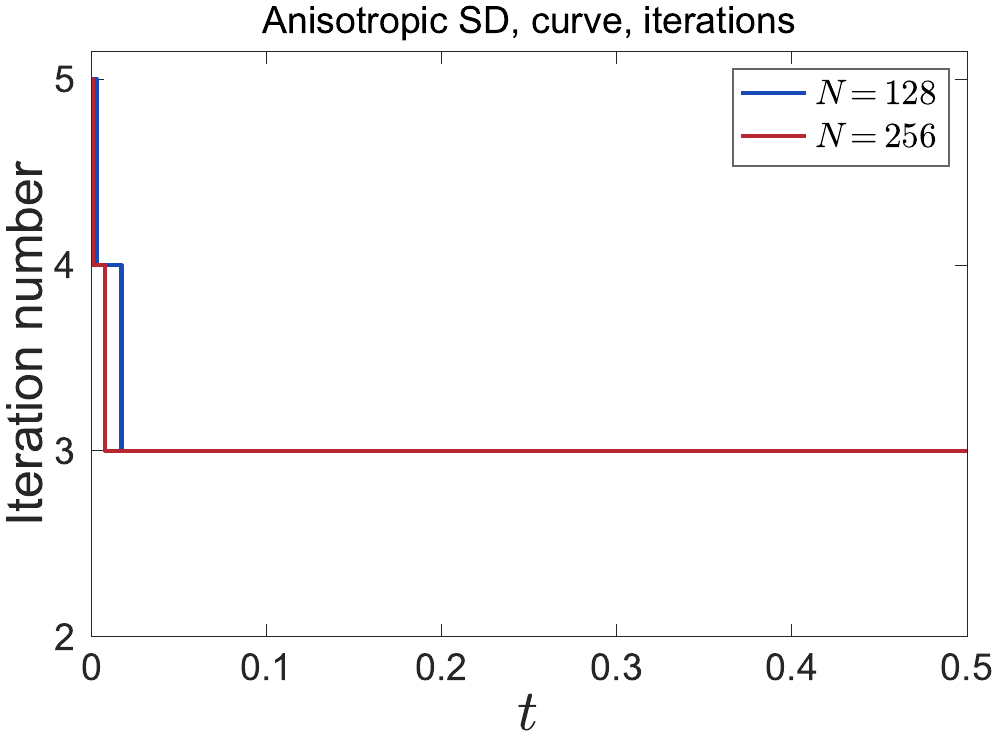}\qquad
\includegraphics[width=0.4\textwidth]{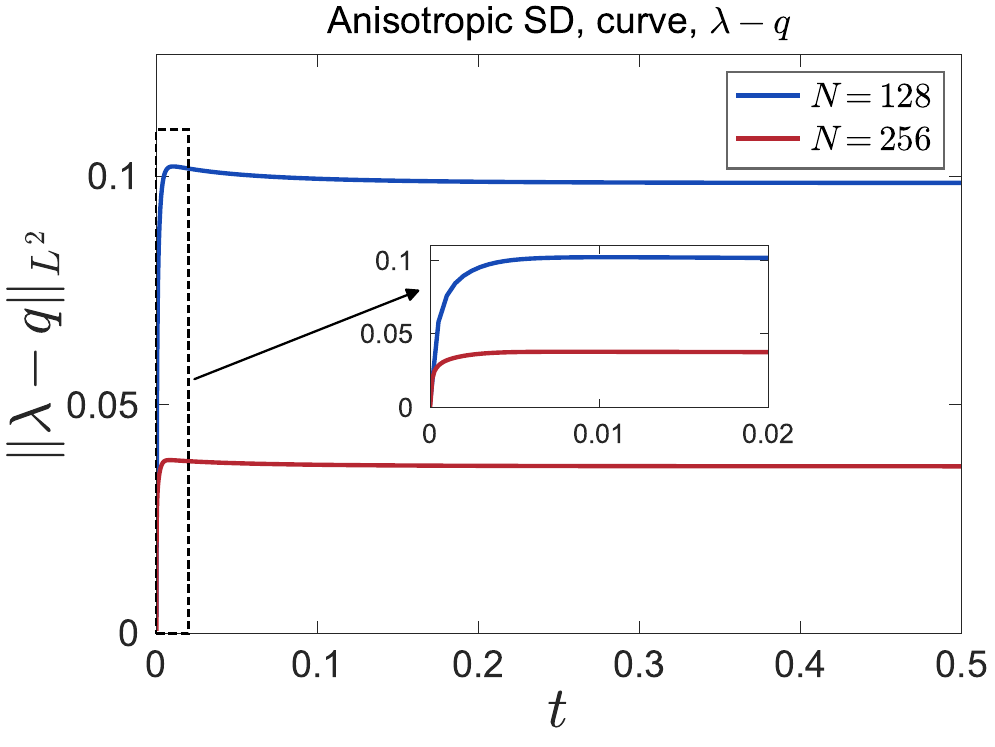}
\includegraphics[width=0.4\textwidth]{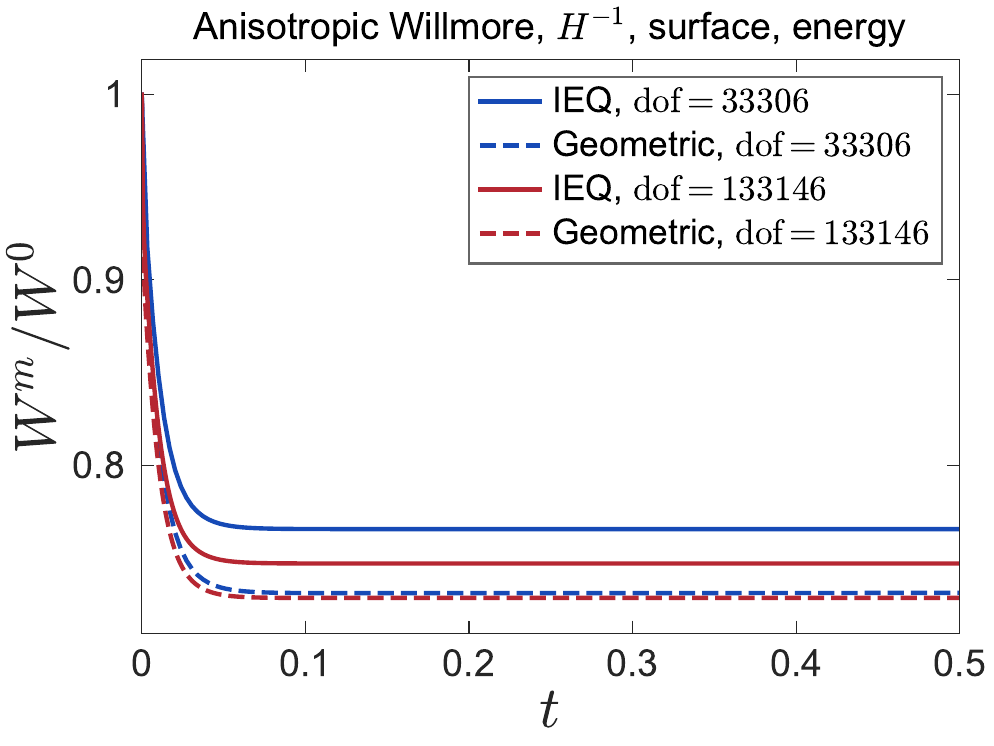}\qquad
\includegraphics[width=0.4\textwidth]{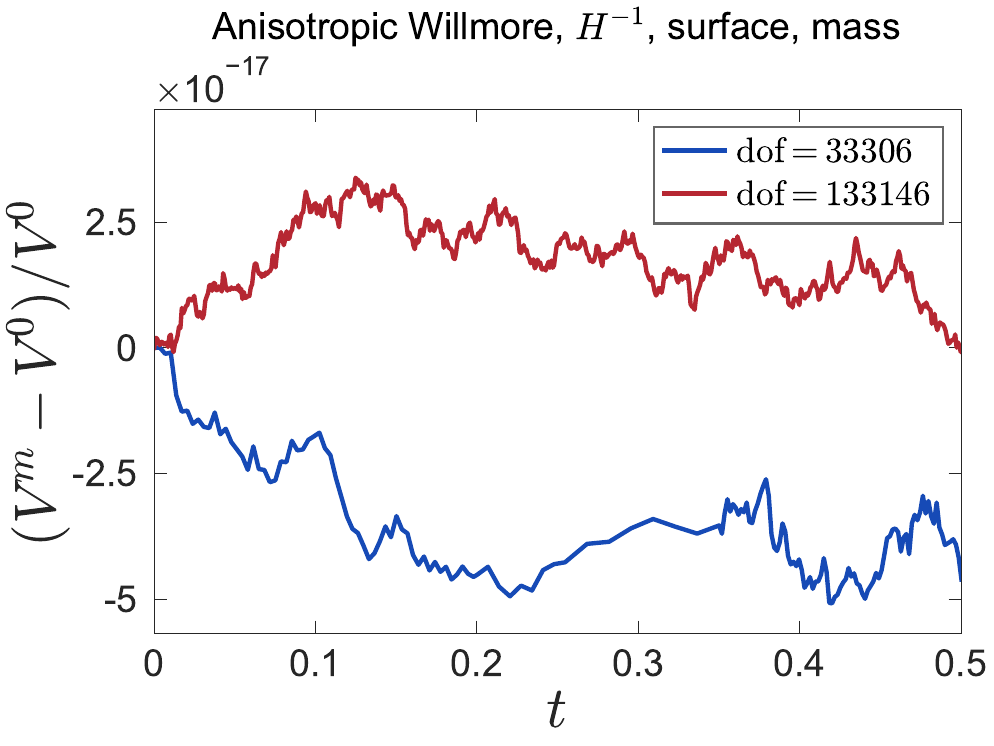}
\includegraphics[width=0.4\textwidth]{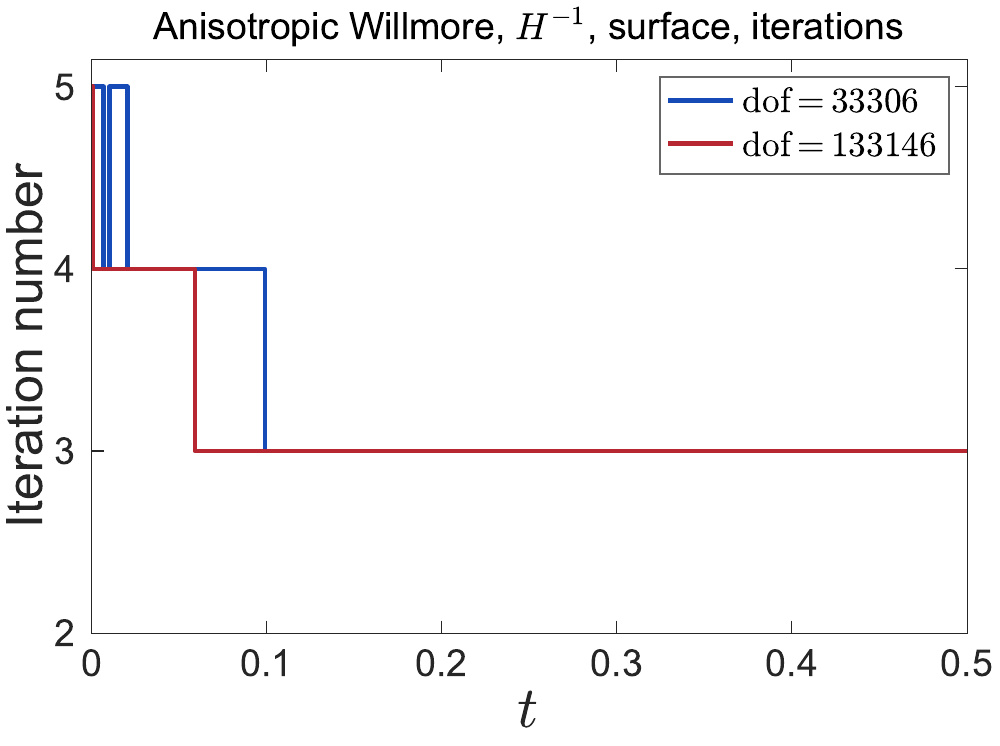}\qquad
\includegraphics[width=0.4\textwidth]{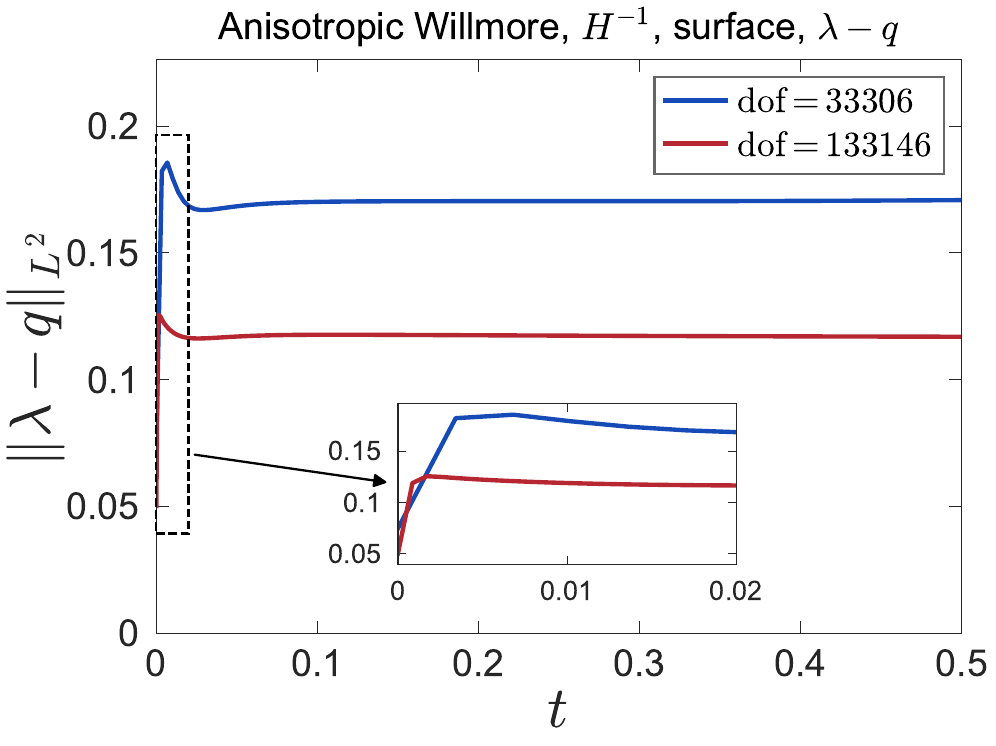}
\caption{The first four correspond to the
anisotropic SD of the flower curve, and the
last four to the $H^{-1}$ anisotropic Willmore flow of Clarenz.  In each
case, the upper row shows the normalized IEQ and geometric energies (left) and
the normalized area change $(A^m-A^0)/A^0$ or volume change $(V^m-V^0)/V^0$
(right), and the lower row shows the number of Newton iterations (left) and
the IEQ defect \eqref{eq:num-defect} (right).}
\label{fig:num-ind-h-1}
\end{figure}

The four indicators are shown in Fig.~\ref{fig:num-ind-l2} for the
$L^2$-gradient flow and in Fig.~\ref{fig:num-ind-h-1} for the
$H^{-1}$-gradient flows, for $0\leq t\leq0.5$.  From these figures, we can
conclude that:

(i) The IEQ energy is monotonically decreasing in all cases, which agrees
with the unconditional energy stability in Theorems~\ref{thm:discrete-energy}
and \ref{thm:h-1-discrete-energy}. The maximal relative difference between the two
energies is $0.92\%$ and $0.24\%$ for the anisotropic MCF, $3.4\%$ and $1.3\%$
for the anisotropic SD, and $4.9\%$ and $2.6\%$ for the anisotropic Willmore
flow, on the coarser and the finer mesh, respectively. As dof increases, the IEQ energy gets closer to the geometric energy. Moreover, the IEQ energy approximates the geometric energy better for the $L^2$ gradient flow.  

(ii) The normalized area loss of the anisotropic SD and the anisotropic Willmore $H^{-1}$-gradient flow is at the level of $10^{-16}$, i.e.\ at the level of machine errors.  This confirms the exact mass conservation in Theorem~\ref{thm:h-1-discrete-energy}.  For
the anisotropic MCF, the area enclosed by the evolving simple closed curve decreases linearly with rate $2\pi$, which preserves the area decay rate. In fact, for $d=2$, $L=1$ and $\beta=\gamma$ in \eqref{eq:num-gamma-G}, one has \cite{MR1248680}: 
$$\frac{\rm d}{{\rm d}t}A(t)=\frac12\int_{\Gt}\gamma(\vv\nu)\varkappa_\gamma\rd\Gamma=\frac{1}{2}(-2{\rm vol}(\partial\mathcal{W}))=-\pi\sqrt{{\rm det}(\mbox{$G_1$})}=-2\pi,$$
for every simple closed curve. 

(iii) Although the schemes are fully implicit, most time steps require three
Newton iterations, and at most five iterations are needed in all tests. Therefore, our numerical methods are weakly nonlinear and computationally efficient.

(iv) The IEQ defect \eqref{eq:num-defect} remains bounded throughout
the evolution and decreases under mesh refinement, with smaller defect observed for the $L^2$-gradient flow.

\FloatBarrier
\subsection{Morphological evolutions}
\label{subsec:num-morphology}

Finally, we consider morphological evolutions of closed curves and surfaces.
We first illustrate the $H^{-1}$-gradient flow of  the Willmore energy \begin{equation*}
  E(\Gamma)=\frac12\int_\Gamma(\varkappa-\overline{\varkappa})^2\rd\Gamma
  =\frac12\int_\Gamma|\vv\nu|(\varkappa-\overline{\varkappa})^2\rd\Gamma,
\end{equation*}
with spontaneous curvature $\overline{\varkappa}=-2$, which conserves the total enclosed area/volume.

For both curves and surfaces, we use the time step $\tau=10^{-5}$. 

\begin{figure}[htp!]
\centering
\setlength{\tabcolsep}{2pt}
\begin{tabular}{@{}ccc@{}}
\includegraphics{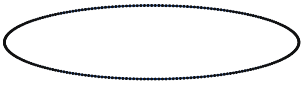} &
\includegraphics{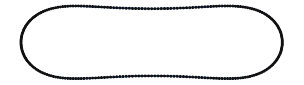} &
\includegraphics{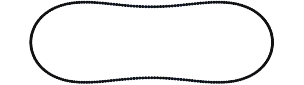} \\
(a) $t=0$ & (b) $t=0.00199$ & (c) $t=0.0102$ \\[6pt]
\includegraphics{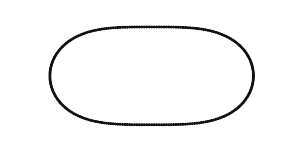} &
\includegraphics{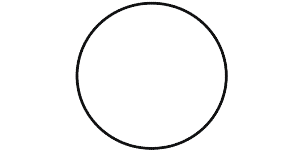} &
\includegraphics{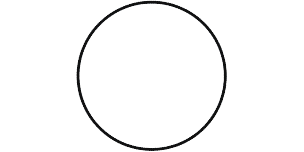} \\
(d) $t=0.0300$ & (e) $t=0.0960$ & (f) $t=0.437$
\end{tabular}
\caption{Morphological evolution of a $2\times0.5$ ellipse under the
$H^{-1}$-gradient flow with spontaneous curvature $\overline{\varkappa}=-2$, $\tau=10^{-5}$.}
\label{fig:num-morph-spontaneous-curve}
\end{figure}

\textbf{$H^{-1}$ Willmore flow: curve.} The initial curve is a
$2\times0.5$ ellipse, with $N=256$ vertices.
Figure~\ref{fig:num-morph-spontaneous-curve} shows that concave parts first develop along the flat sides with small curvature. It subsequently becomes rounder and approaches a circle.

\begin{figure}[htp!]
\centering
\setlength{\tabcolsep}{2pt}
\begin{tabular}{@{}ccc@{}}
\includegraphics{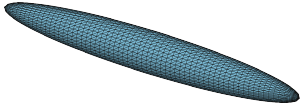} &
\includegraphics{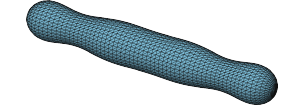} &
\includegraphics{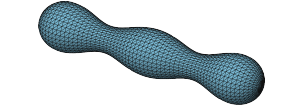} \\
(a) $t=0$ & (b) $t=0.00275$ & (c) $t=0.0244$ \\[6pt]
\includegraphics{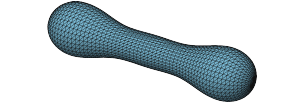} &
\includegraphics{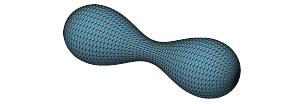} &
\includegraphics{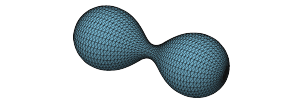} \\
(d) $t=0.102$ & (e) $t=0.199$ & (f) $t=1.95$
\end{tabular}
\caption{Morphological evolution of a $4\times0.5\times0.5$ ellipsoid under
the $H^{-1}$-gradient flow with spontaneous curvature
$\overline{\varkappa}=-2$, up to $t=1.95$, with $\tau=10^{-5}$.}
\label{fig:num-morph-spontaneous-surface}
\end{figure}

\textbf{$H^{-1}$ Willmore flow: surface.} The initial surface is a
$4\times0.5\times0.5$ ellipsoid with $N=2018$ vertices. As shown in
Fig.~\ref{fig:num-morph-spontaneous-surface}, it first develops two concave parts along the flat sides. Then the two concave parts merge to a neck, and it approaches a nonspherical stationary shape with two spheres joined by a thin neck.

\FloatBarrier

To illustrate the desirable behaviour of our numerical schemes, we now present two challenging situations for surfaces: a strongly anisotropic energy, whose energy density is not convex and is regularized by a mean curvature term; and anisotropic Willmore flows with a nearly crystalline BGN anisotropy. In these two examples, the initial surface is discretized by the icosahedral meshes of Section~\ref{subsec:num-convergence} with $N=2562$, and the gIEQ--MDR schemes are used.

\textbf{Strong anisotropy with regularization.}  Following \cite{Zhang26coupled}, we consider the energy
\begin{equation}
  E(\Gamma)=\int_\Gamma\Bigl[\gamma_c(\vv\nu)+\frac{\varepsilon^2}{2}\varkappa^2\Bigr]\rd\Gamma,
  \qquad
  \gamma_c(\vv p)=|\vv p|\Bigl(1+a\sum_{i=1}^3\Bigl(\frac{p_i}{|\vv p|}\Bigr)^3\Bigr),
  \label{eq:num-zhang}
\end{equation}
with $a=\frac12$ and $\varepsilon^2=0.05$. The cubic anisotropy $\gamma_c$ is convex if and
only if $|a|\leq \sqrt6/8\approx0.306$.  Hence $a=\frac12$ is strongly anisotropic, which
is regularized by the mean curvature. Here we have $f(\varkappa,\vv p)=2\gamma_c(\vv p)+\varepsilon^2|\vv
p|\varkappa^2$ and $q=\sqrt f$.   

We consider the $L^2$-gradient flows and the
$H^{-1}$-gradient flows with $\beta=1$, i.e., the strongly
anisotropic MCF and SD, respectively.  The initial surface is a $1.4 \times 1 \times 0.8$ ellipsoid \cite{Zhang26coupled}, with $\tau=3.81\times10^{-5}$.

\begin{figure}[htp!]
\centering
\setlength{\tabcolsep}{2pt}
\begin{tabular}{@{}ccc@{}}
\includegraphics[width=0.3\linewidth]{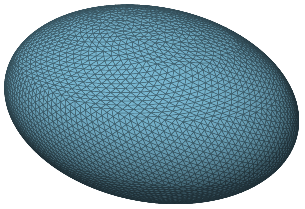} &
\includegraphics[width=0.3\linewidth]{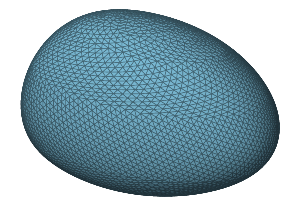} &
\includegraphics[width=0.3\linewidth]{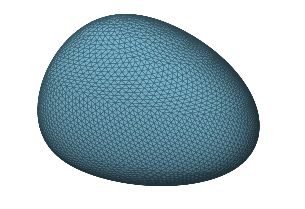} \\
(a) $t=0$ & (b) $t=0.05$ & (c) $t=0.1$ \\[6pt]
\includegraphics[width=0.3\linewidth]{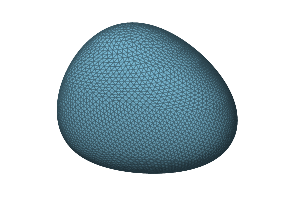} &
\includegraphics[width=0.3\linewidth]{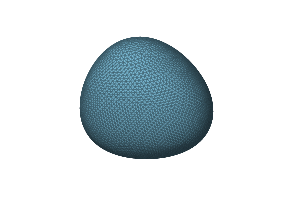} &
\includegraphics[width=0.3\linewidth]{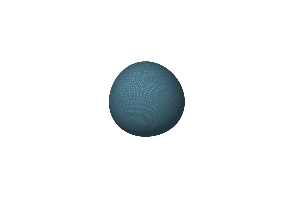} \\
(d) $t=0.15$ & (e) $t=0.2$ & (f) $t=0.25$
\end{tabular}
\caption{Morphological evolution of a $1.4 \times 1 \times 0.8$ ellipsoid under the strongly anisotropic MCF, with $\tau=3.81\times10^{-5}$.}
\label{fig:num-morph-strong-mcf}
\end{figure}

\begin{figure}[htp!]
\centering
\setlength{\tabcolsep}{2pt}
\begin{tabular}{ccc}
\includegraphics[width=0.3\linewidth]{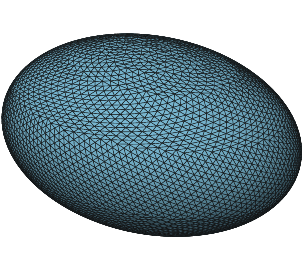} &
\includegraphics[width=0.3\linewidth]{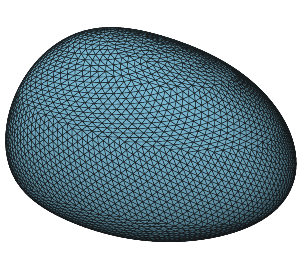} &
\includegraphics[width=0.3\linewidth]{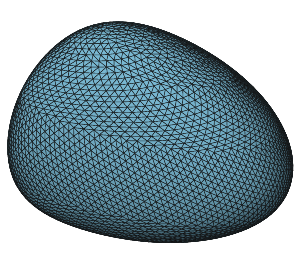} \\
(a) $t=0$ & (b) $t=7.51\times10^{-3}$ & (c) $t=1.50\times10^{-2}$ \\[6pt]
\includegraphics[width=0.3\linewidth]{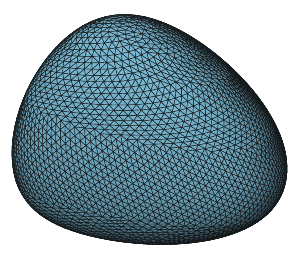} &
\includegraphics[width=0.3\linewidth]{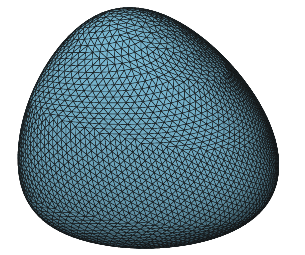} &
\includegraphics[width=0.3\linewidth]{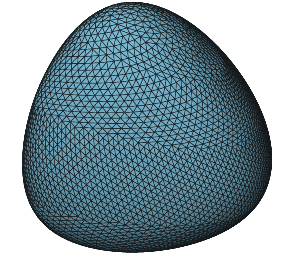} \\
(d) $t=3.00\times10^{-2}$ & (e) $t=5.99\times10^{-2}$ & (f) $t=0.25$
\end{tabular}
\caption{Morphological evolution of a $1.4 \times 1 \times 0.8$ ellipsoid under the strongly anisotropic SD, with $\tau=3.81\times10^{-5}$.}
\label{fig:num-morph-strong-sd}
\end{figure}

The morphological evolutions are shown in
Figs.~\ref{fig:num-morph-strong-mcf} and \ref{fig:num-morph-strong-sd}.  We
observe that the evolutions are not symmetric with respect to the origin,
since the anisotropy is not symmetric, i.e.\ $\gamma_c(-\vv p)\neq\gamma_c(\vv p)$. For the strongly anisotropic MCF, the surface shrinks to a tetrahedron-like shape, but it becomes rounder as the mean curvature term in \eqref{eq:num-zhang} becomes
dominant.  For the strongly anisotropic SD, the surface evolves to a tetrahedron-like shape at the equilibrium, where the anisotropic energy is the dominant effect. The mesh points are well-behaved in both simulations, and no mesh regularization is needed.

\textbf{Crystalline BGN anisotropy.}  We consider the anisotropic Willmore flow with the nearly crystalline BGN anisotropy \eqref{eq:num-gamma-G} with $L=3$ and
$G_\ell=\epsilon^2I+(1-\epsilon^2)\,\vv e_\ell\otimes\vv e_\ell$
\cite{BGN08ani3d}, where $I$ is the identity matrix and $\vv e_\ell$ is the
$\ell$-th unit vector, i.e. \cite{Barrett20}
\begin{equation}
  \gamma_\epsilon(\vv p)=\sum_{\ell=1}^3\sqrt{\epsilon^2|\vv p|^2+(1-\epsilon^2)p_\ell^2},
  \qquad\epsilon=0.1,
  \label{eq:num-sum-three}
\end{equation}
which is a smooth approximation of $|p_1|+|p_2|+|p_3|$.  Its Wulff shape is
a cube with rounded edges and corners.  We compare the energies of Clarenz,
$\alpha(\vv p)=|\vv p|$ and $\beta=1$, and of Pozzi, $\alpha=\beta=\gamma_\epsilon$, in
\eqref{eq:anisotropic-energy} and \eqref{eq:anisotropic-normal-velocity}.  The
initial surface is the unit sphere, with $\tau=9.53\times10^{-6}$.

\begin{figure}[htp!]
\centering
\begin{tabular}{@{}c@{\hspace{2pt}}c@{\hspace{0.03\textwidth}}c@{\hspace{2pt}}c@{}}
Clarenz & Pozzi & Clarenz & Pozzi \\[2pt]
\includegraphics[width=0.2\linewidth]{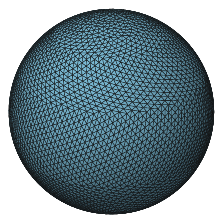} &
\includegraphics[width=0.2\linewidth]{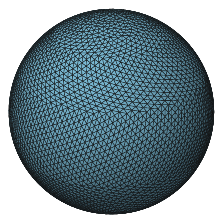} &
\includegraphics[width=0.2\linewidth]{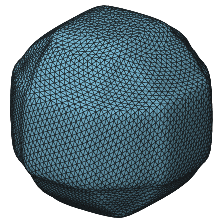} &
\includegraphics[width=0.2\linewidth]{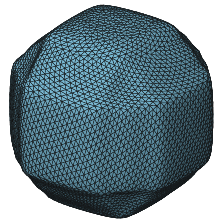} \\
\multicolumn{2}{c@{\hspace{0.03\textwidth}}}{(a) $t=0$} &
\multicolumn{2}{c}{(b) $t=2.86\times10^{-4}$} \\[6pt]
\includegraphics[width=0.2\linewidth]{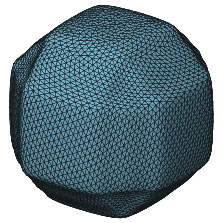} &
\includegraphics[width=0.2\linewidth]{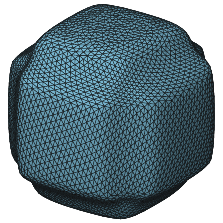} &
\includegraphics[width=0.2\linewidth]{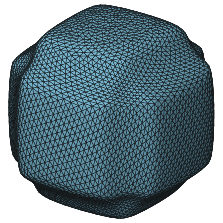} &
\includegraphics[width=0.2\linewidth]{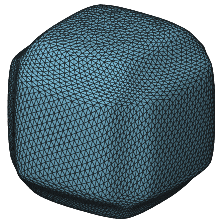} \\
\multicolumn{2}{c@{\hspace{0.03\textwidth}}}{(c) $t=6.96\times10^{-4}$} &
\multicolumn{2}{c}{(d) $t=1.51\times10^{-3}$} \\[6pt]
\includegraphics[width=0.2\linewidth]{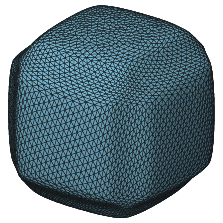} &
\includegraphics[width=0.2\linewidth]{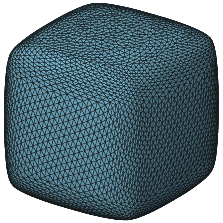} &
\includegraphics[width=0.2\linewidth]{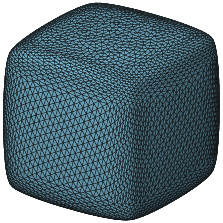} &
\includegraphics[width=0.2\linewidth]{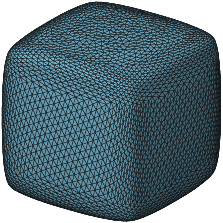} \\
\multicolumn{2}{c@{\hspace{0.03\textwidth}}}{(e) $t=3.48\times10^{-3}$} &
\multicolumn{2}{c}{(f) $t=0.02$} 
\end{tabular}
\caption{Morphological evolutions of the unit sphere under the anisotropic
Willmore flows of Clarenz and Pozzi with $\gamma_\epsilon$ in
\eqref{eq:num-sum-three}, $\epsilon=0.1$, and
$\tau=9.53\times10^{-6}$.  Each pair shows the two flows at the same time.}
\label{fig:num-morph-willmore}
\end{figure}

Fig.~\ref{fig:num-morph-willmore} compares the two evolutions at the same
times.  In both cases, flat faces and sharp corners appear on the sphere during the
evolution, and the two flows evolve to the same equilibrium shape (cf.\ Fig.~\ref{fig:num-morph-willmore}~(f)), which agrees with the theoretical result that both energies attain their global minimum at Wulff shapes
\cite[Theorem~4.1]{Clarenz04wulff}, \cite[Theorem~2.4]{Pozzi15anisoWillmore}.
However, the two flows evolve with different speeds. The anisotropic Willmore flow of Pozzi is faster, and it has already reached the equilibrium in
Fig.~\ref{fig:num-morph-willmore}~(e), while the flow of Clarenz has not.
Moreover, the equilibrium shape of the two flows is stable.  As before, the mesh points are
well-behaved, and no mesh regularization is required.

\section{Conclusions}
\label{sec:con}

We have developed a geometric-IEQ framework for weighted $L^2$-gradient flows of a class of curvature-dependent energies of closed curves and surfaces. The framework separates the energy dissipation in the normal direction from tangential mesh motion. Coupling it with the BGN or MDR formulation yields parametric finite element schemes with unconditional decay of the discrete quadratic IEQ energy. The framework also extends to anisotropic Willmore flow and volume conserving $H^{-1}$-gradient flows, with exact area or volume conservation for the latter.

Numerical experiments support the theoretical results and demonstrate the effectiveness of the proposed methods in capturing complex geometric evolutions while maintaining good mesh quality.

\section*{Acknowledgement}

This work was partially supported by the National Natural Science Foundation of China
(No. 12401572, Q.Z.).

\bibliographystyle{model1b-num-names} 
\bibliography{bib}

\end{document}